\documentclass[11pt]{article}
\usepackage[latin9]{inputenc}
\usepackage{amsmath}
\usepackage{amssymb}

\makeatletter
\usepackage[latin9]{luainputenc}
\usepackage{xcolor}
\usepackage{pdfcolmk}
\usepackage{float}
\usepackage{mathrsfs}
\usepackage{amsmath}
\usepackage{amssymb}
\usepackage{graphicx}
\usepackage{esint}
\usepackage{authblk}
\PassOptionsToPackage{normalem}{ulem}
\usepackage{ulem}

\makeatletter

\providecolor{lyxadded}{rgb}{0,0,1}
\providecolor{lyxdeleted}{rgb}{1,0,0}
\usepackage{amsfonts}
\usepackage{graphicx}
\usepackage{caption}
\usepackage{subcaption}
\usepackage{amsfonts}
\usepackage{amsthm}

\theoremstyle{definition}
\newtheorem{theorem}{Theorem} 
\newtheorem{corollary}{Corollary} 
\newtheorem{lemma}{Lemma}
\newtheorem{proposition}{Proposition}
 
\theoremstyle{remark}

\newtheorem{remark}{Remark}

\makeatother

\begin{document}
\title{The Euler Scheme for Fractional Stochastic Delay Differential Equations
	with Additive Noise}
\author{Orimar Sauri \thanks{osauri@math.aau.dk}}
\affil{Department of Mathematical Sciences, Aalborg University, Denmark}
\maketitle
	\begin{abstract}
		In this paper we consider the Euler-Maruyama scheme for a class of
		stochastic delay differential equations driven by a fractional Brownian
		motion with index $H\in(0,1)$. We establish the consistency of the
		scheme and study the rate of convergence of the normalized error process.
		This is done by checking that the generic rate of convergence of the
		error process with stepsize $\Delta_{n}$ is $\Delta_{n}^{\min\{H+\frac{1}{2},3H,1\}}$.
		It turned out that such a rate is suboptimal when the
		delay is smooth and $H>1/2$. In this context, and in contrast to
		the non-delayed framework, we show that a convergence of order $H+1/2$
		is achievable.
	\end{abstract}
	
	% Use the \maketitle command after the abstract

	%% Beginning of text
	
	\section{Introduction\label{intro}}
	
	\subsection*{Overview}
	
	This paper is concerned with numerical approximations of solutions
	of stochastic delay differential equations (SDDE from now on) of the
	form 
	\begin{equation}
		X_{t}=\begin{cases}
			x_{0}(0)+\int_{0}^{t}\int_{[0,\tau]}b(X_{s-r})\eta(\mathrm{d}r)\mathrm{d}s+B_{t}, & t\geq0,\\
			x_{0}(t), & -\tau\leq t<0.
		\end{cases}\label{SDDE-1}
	\end{equation}
	Above $\eta$ represents a finite signed measure on $[0,\tau]$, $\tau>0$,
	and $B$ denotes a fractional Brownian motion with index $H\in(0,1)$
	(for completeness we also consider the case $H=1/2$). We are particularly
	interested in the optimal rate of convergence of the Euler scheme
	associated to (\ref{SDDE-1}). To obtain such a rate, in this work
	we focus on deriving non-trivial limit theorems for the (normalized)
	error process linked to this method.
	
	SDDEs are often seen as the natural generalization of classic SDEs
	to the non-Markovian framework which, in the context of (\ref{SDDE-1}),
	is carry through $\eta$. In this regard, the more weight $\eta$
	assigns to values near to $\tau$, the more the influence the past
	has over the current state of the process. This simple way of introducing
	memory, or \textit{long-range dependence}, makes SDDEs potential candidates
	to model relevant financial quantities such as interest rates and
	stochastic volatility. See for instance \cite{Rihan2021,KuchlerPlaten07,BaoYinYuan16,YuriSwWu05,YuriSwWu07,ArrioHuMohPap07}
	and references therein. Another distinctive characteristic of solutions
	of (\ref{SDDE-1}), which is the main motivation of this work, is
	that they allow for \textit{roughness} in the sense of \cite{GathJaiRos18}
	(see also \cite{BenLundPakk16}), making them suitable for modelling
	\textit{roughness and long memory on volatility.} For a survey on
	the modeling of stochastic volatility we refer the reader to \cite{DiNunnoKubMishuYour23}.
	See also \cite{BayerfrizFukaGatherakJacqRosem23}. 
	In the context of option pricing 
	\subsection*{Related work}
	
	The literature on the limiting behavior of the error of numerical
	methods for classical SDEs (without delay) driven by semimartingales
	is vast and goes back to the seminal works of \cite{KurtzPotter91}
	and \cite{JacPrott98}. It is remarkable that this is also the case
	when the SDE is driven by a fractional Brownian motion (fBm) rather
	than a semimartingale. We would like to emphasize that most of the
	existing works in this framework concentrate on the non-rough set-up,
	i.e. the situation in which the Hurst exponent of the driving fBm
	is above $1/2$. See for instance \cite{Neuenkirch06,HuLiuNualart16,ZhangYuan21}.
	For papers dealing with the rough case we refer the reader to \cite{NeueNourdin07,GradNourdin09,LiuTindel19}.
	
	In contrast, when a delay is added, references are more scarce in
	both frameworks. Current research seems to focus on the rate of convergence
	of the scheme rather than in the limit distribution of the error process.
	For the semimartingale set-up see for instance \cite{BUCKWAR00,HuMohaYan04,HoffmanMuller06,ZhangGandHu09,KloedenShardlow12,BaoYinYuan16}
	and for numerical schemes for fractional SDDEs we refer to \cite{GarzonTindelTorres19,MAHMOUDITahmaseb22}.
	Despite the above, to the best of our knowledge, the problem considered
	in this work has not been addressed anywhere else although it can
	be seen as an extension of the works in \cite{HuLiuNualart16,ZhangYuan21}
	and \cite{LiuTindel19}.
	
	\subsection*{Main contributions}
	
	This paper establishes strong rates of convergence as well as limit
	theorems for the error process resulting from applying the Euler-Maruyama
	method to (\ref{SDDE-1}). Specifically, under a smooth and linear growth condition on the drift component $b$, we show
	that the \textit{generic rate} of convergence is $\Delta_{n}^{\min\{H+\frac{1}{2},3H,1\}}$,
	that is, for every $T>0$ and $p\geq1$ there is a constant $C>0$
	independent of $n\in\mathbb{N}$, such that
	\[
	\sup_{0\leq t\leq T}\mathbb{E}(\vert U_{t}^{n}\vert^{p})^{~1/p}\leq C\Delta_{n}^{\min\{H+\frac{1}{2},3H,1\}},\,\,\,H\in(0,1),
	\]
	where 
	\[
	U_{t}^{n}=X_{t}-X_{t}^{n},
	\]
	in which $X$ is the solution of (\ref{SDDE-1}) and $X^{n}$ its
	the Euler scheme (see (\ref{eq:Eulermarscheme}) below for a detailed
	definition). 
	Additionally, when $H\geq1/2$, we prove the following: 
	\begin{enumerate}
		\item There is a bias process $\mathfrak{B}^{n}$, such that $(U^{n}-\mathfrak{B}^{n})/\Delta_{n}$
		converges uniformly in compacts towards a non-degenerate stochastic process. However, in general, the sequence
		$\mathfrak{B}^{n}/\Delta_n$ is not convergent.
		\item If $\eta$ admits a continuous density $U^{n}/\Delta_n$ is asymptotically negligible. Furthermore, in this situation,  a rate of convergence of
		order $H+1/2$ can be attained.
	\end{enumerate}
	\subsection*{Structure of the paper}
	
	The organization of this work is as follows. Section 2 focus on introducing
	the main mathematical concepts and some basic results that will be
	used through the paper. We also present some aspects of Malliavin
	calculus due to its relevance in our proofs. In Section 3 we introduce
	the Euler-Maruyama scheme associated to (\ref{SDDE-1}) and present
	our main results for the error process. For the sake of exposition,
	we postpone all our proofs to the end of the paper, that is, to Section
	4. 
	
	\section{Preliminaries and basic results \label{sec:Preliminaries-and-basic}}
	
	This part is devoted to introduce our set-up as well as state some
	basic results that will be used later.
	
	\subsection{Basic convergence concepts and fractional Brownian motions \label{subsec:Linear-fractional--stable}}
	
	As it is costumary, $(\Omega,\mathcal{F},\mathbb{P})$ will represent
	a complete probability space. The symbols $\overset{\mathbb{P}}{\rightarrow}$
	and $\overset{d}{\rightarrow}$ stand, respectively, for convergence
	in probability and distribution of random vectors (r.v.'s for short).
	For a sequence of random vectors $(\xi_{n})_{n\geq1}$ defined on
	$\left(\Omega,\mathcal{F},\mathbb{P}\right)$, we set, respectively,
	$\xi_{n}=\mathrm{o}_{\mathbb{P}}(1)$ and $\xi_{n}=\mathrm{O}_{\mathbb{P}}(1)$
	whenever $\xi_{n}\overset{\mathbb{P}}{\rightarrow}0$ or if $\xi_{n}$
	is bounded in probability, respectively. Let $(H_{t}^{n})_{t\geq0,n\in\mathbb{N}}$
	be a sequence of càdlàg processes defined on $\left(\Omega,\mathcal{F},\mathbb{P}\right)$.
	We will write $H^{n}\overset{u.c.p}{\longrightarrow}H$ if $H^{n}$
	converges uniformly on compacts in probability towards $H$. Now,
	given a sub-$\sigma$-field $\mathcal{G}\subseteq\mathcal{F}$ and
	a random vector $\xi$ (defined possibly on an extension of $\left(\Omega,\mathcal{F},\mathbb{P}\right)$)
	we say $\xi_{n}$ converges\textit{ $\mathcal{G}$-stably in distribution}
	towards $\xi$, in symbols $\xi_{n}\overset{\mathcal{G}\text{-}d}{\longrightarrow}\xi$,
	if for any $\mathcal{G}$-measurable random variable $\zeta$, $(\xi_{n},\zeta)\overset{d}{\rightarrow}(\xi,\zeta)$,
	as $n\rightarrow\infty$. In this framework, if $(H_{t}^{n})_{t\in\mathscr{T},n\in\mathbb{N}}$,
	$\mathscr{T}\subseteq\mathbb{R}$, is a family of stochastic processes,
	we will write $H^{n}\overset{\mathcal{G}\text{-}fd}{\longrightarrow}H$
	if the finite-dimensional distributions (f.d.d. for short) of $H^{n}$
	converge $\mathcal{G}$-stably toward the f.d.d. of $X$. Furthermore,
	if $(H_{t}^{n})_{0\leq t\leq T,n\in\mathbb{N}}$ is a sequence of
	continuous processes, we write $H^{n}\overset{\mathcal{G}\text{-}\mathcal{C}([0,T])}{\Longrightarrow}H$,
	if $H^{n}$ converges weakly to $H$ in the uniform topology and $H^{n}\overset{\mathcal{G}\text{-}fd}{\longrightarrow}H$.
	We refer the reader to \cite{HauslerLuschgy15} for a concise exposition
	of stable convergence.
	
	In this work $(B_{t})_{t\in\mathbb{R}}$ will denote a \textit{fractional
		fractional Brownian motion} (fBm from now on), i.e. if its covariance
	function can be represented as 
	\begin{equation}
		B_{t}=\int_{-\infty}^{t}[(t-s)_{+}^{H-1/2}-(-s)_{+}^{H-1/2}]\mathrm{d}W_{s}\text{, \ \ }t\in\mathbb{R}\text{,}\label{fBmdef}
	\end{equation}
	where $(x)_{+}^{p}:=x^{p}\mathbf{1}_{x>0}$, $H\in(0,1)$ and $W$
	is a two-sided Wiener process defined on $(\Omega,\mathcal{F},\mathbb{P})$.
	Note that we are also considering the case $H=1/2$, i.e. we are also
	consider the situation in which $B=W$. It is well-known that the
	process $Z$ has $\lambda$-Hölder continuous paths for any $\lambda<H$.
	More precisely, for every $T>0$ and $0<\lambda<H$, there is a positive
	random variable, say $\xi_{\lambda,T}$ , such that 
	\begin{equation}
		\lvert B_{t}-B_{s}\lvert\leq\xi_{\lambda,T}\lvert t-s\lvert^{\lambda},\,\,\forall\,t,s\in[0,T].\label{eq:HolderexplicitfBm}
	\end{equation}
	Furthermore, $\xi_{\lambda,T}$ has finite moments of all orders and
	is independent of $t,s\in[0,T]$. For more details we refer the reader
	to \cite{AzmooSotViitYa14} and references therein. Let $\beta=H-1/2$.
	Note that for $t\geq0$ we can decompose $Z$ as 
	\begin{equation}
		B_{t}=V_{t}+Z_{t},\label{eq:decompFLSM}
	\end{equation}
	where $V_{0}=Z_{0}=0$ and for $t>0$ 
	\[
	Z_{t}=\int_{0}^{t}(t-s)^{\beta}\mathrm{d}W_{s},\,\,\,V_{t}=\int_{(-\infty,0]}\left[(t-s)^{\beta}-(-s)^{\beta}\right]\mathrm{d}W_{s}.
	\]
	Since 
	\[
	\int_{0}^{t}\left(\int_{(-\infty,0]}(u-s)^{(\beta-1)2}\mathrm{d}s\right)^{1/2}\mathrm{d}u<\infty,\,\,t>0,
	\]
	then by the Stochastic Fubini Theorem (see e.g. \cite{BasseBN11}),
	the process $(V_{t})_{t\geq0}$ admits an absolutely continuous version
	given by 
	\begin{equation}
		V_{t}=\int_{0}^{t}v_{u}\mathrm{d}u,\,\,\,v_{u}:=\int_{(-\infty,0]}(u-s)^{(\beta-1)}\mathrm{d}W_{s}.\label{eq:acpartZ}
	\end{equation}

	\subsection{Differential resolvents\label{subsec:Differential-resolvents}}
	
	The proof of our main results rely on differential resolvents of the
	so-called Volterra measure kernels. Therefore, in this part we discuss
	such concepts and their basic properties. Fix a measurable and bounded
	process $(b_{t})_{t\geq0}$ and let 
	\[
	\kappa(t,A,\omega)=-\int_{[0,t]}b_{t-r}(\omega)\mathbf{1}_{A}(t-r)\eta(\mathrm{d}r),\,\,t\geq0,\,\,A\in\mathcal{B}([0,+\infty)),
	\]
	where $\eta$ is a finite (deterministic) signed measure with support
	on $[0,\tau]$, $0\leq\tau<\infty$. Following the terminology of
	\cite{GripenbergLondenStaffans90}, Chapter 10, for every $\omega\in\Omega$,
	$\kappa(\cdot,\omega)$ is a \textit{Volterra kernel measure} of type
	$B^{\infty}$ on $[0,+\infty)$, i.e. $\forall\,t\geq0$, $\kappa(t,\cdot,\omega)$
	is a signed finite measure with support on $[0,t]$, for every Borelian
	set $A$, $\kappa(\cdot,A,\omega)$ is measurable, and the mapping
	$t\mapsto\lvert\kappa\lvert(t,[0,+\infty),\omega)$ is bounded. Therefore,
	for every $\omega\in\Omega$ there is a locally bounded and measurable
	mapping $R_{\omega}^{\kappa}:[0,+\infty)^{2}\rightarrow\mathbb{R}$,
	known as the \textit{differential resolvent} of $\kappa(\cdot,\omega)$,
	that satisfies the following: For every $t,s\geq0$, $R_{\omega}^{\kappa}(t,s)=0$
	if $s>t$, $R_{\omega}^{\kappa}(t,t)=1$, otherwise 
	\begin{equation}
		\begin{aligned}R_{\omega}^{\kappa}(t,s) & =\mathbf{1}_{0\leq s\leq t}+\int_{0}^{t}\int_{[0,u]}b_{u-r}(\omega)R_{\omega}^{\kappa}(u-r,s)\eta(\mathrm{d}r)\mathrm{d}u\\
			& =\mathbf{1}_{0\leq s\leq t}+\int_{0}^{t}R_{\omega}^{\kappa}(t,u)\int_{[0,u]}b_{u-r}(\omega)\mathbf{1}_{s\leq u-r}\eta(\mathrm{d}r)\mathrm{d}u.
		\end{aligned}
		\label{eq:resolventequation}
	\end{equation}
	Furthermore, the mapping $t\in[s,+\infty)\mapsto R_{\omega}^{\kappa}(t,s)$
	($s\mapsto R_{\omega}^{\kappa}(t,s)$) is locally absolutely continuous
	(left-continuous and locally of bounded variation on $[0,+\infty)$).
	We can say a bit more about $R_{\omega}^{\kappa}$. The first part
	of (\ref{eq:resolventequation}), the boundedness of $b$ and Gronwall's
	inequality gives that for every $T>0$ there is $C>0$ independent
	of $\omega$, such that 
	\begin{equation}
		\sup_{0\leq s\leq t\leq T}\vert R_{\omega}^{\kappa}(t,s)\vert\leq e^{CT}.\label{eq:BoundedResolvente}
	\end{equation}
	Additionally, the second part of equation (\ref{eq:resolventequation})
	gives that for any $0\leq a\leq c\leq t$ 
	\begin{equation}
		-\int_{a}^{b}\int_{[0,t-s]}b_{s}(\omega)R_{\omega}^{\kappa}(t,s+r)\eta(\mathrm{d}r)\mathrm{d}s=R_{\omega}^{\kappa}(t,c)-R_{\omega}^{\kappa}(t,a).\label{eq:diffresolv}
	\end{equation}
	Therefore, the following holds for every $t>0$: 
	\begin{enumerate}
		\item The mapping $s\in[0,t]\mapsto R_{\omega}^{\kappa}(t,s)$ is Lipschitz
		continuous with derivative 
		\begin{equation}
			\frac{\partial}{\partial s}R_{\omega}^{\kappa}(t,s)=-\int_{[0,t-s]}b_{s}(\omega)R_{\omega}^{\kappa}(t,s+r)\eta(\mathrm{d}r).\label{eq:derivative_diif_resolv}
		\end{equation}
		\item For each $t\geq0$, $R_{\omega}^{\kappa}(t,\cdot)$ induces a unique
		signed measure on $[0,t]$ given by 
		\begin{equation}
			R_{\omega}^{\kappa}(t,\mathrm{d}s)=-\int_{[0,t-s]}b_{s}(\omega)R_{\omega}^{\kappa}(t,s+r)\eta(\mathrm{d}r)\mathrm{d}s.\label{eq:measureinducedbyResolvent}
		\end{equation}
	\end{enumerate}
	Be aware that the construction of $R_{\omega}^{\kappa}$ is done $\omega$-by-$\omega$
	and does not guarantee its measurability as function of $\omega$.
	Afortiori, the next result (whose proof is postponed to Section \ref{sec:Proofs})
	shows that the measurability is preserved due to the joint measurability
	of $b$.
	
	\begin{proposition}\label{propmeasurableresolv} For every $T>0$,
		the mapping $(t,s,\omega)\in[0,T]^{2}\times\Omega\mapsto R_{\omega}^{\kappa}(t,s)$
		is $\mathcal{B}([0,T]^{2})\otimes\mathcal{F}\backslash\mathcal{B}(\mathbb{R})$
		measurable. Furthermore, for every fixed $t\geq0$, the process $(R^{\kappa}(t,s))_{0\leq s\leq t}$
		is bounded and Lipschitz continuous.
		
	\end{proposition}
	
	\begin{remark}\label{remarkboundedlocallyresolvent}Note that from
		(\ref{eq:BoundedResolvente}), for every $T>0,$ the mapping $(\omega,t,s)\in\Omega\times[0,T]^{2}\mapsto R_{\omega}^{\kappa}(t,s)$
		is bounded. A combination of the later property and (\ref{eq:diffresolv})
		let us conclude that for a given $T>0$ and every $0\leq a\leq b\leq t\leq T$
		there exists a contant $C$ only depending on $T$ and $\eta$, such
		that 
		\begin{equation}
			\mid R_{\omega}^{\kappa}(t,b)-R_{\omega}^{\kappa}(t,a)\mid\leq C(b-a).\label{eq:semiglobalLipschitzcondresolvent}
		\end{equation}
		
	\end{remark}
	
	\begin{remark} Let $(\mathcal{F}_{t}^{W})_{t\in\mathbb{R}}$ be the
		completion of the natural filtration of the underlying Wiener process
		$W$. Thus, if the process $b$ is continuous and adapted to such
		a filtration, then we also have that the random field $(R^{\kappa}(t,s))_{0\leq s,t\leq T}$
		is $\mathcal{B}([0,T]^{2})\otimes\mathcal{F}_{T}^{W}\backslash\mathcal{B}(\mathbb{R})$
		measurable for all $T>0$.
		
	\end{remark}
	
	\subsection{Elements of Malliavin Calculus}
	
	In this part we recall some basic definitions and facts of Malliavin
	calculus. Let $\mathcal{H}$ be a real separable Hilbert space with inner
	product $\langle,\rangle_{\mathcal{H}}$ and $\{W(h):h\in\mathcal{H}\}$
	an \textit{isonormal Gaussian process} on $(\Omega,\mathcal{F},\mathbb{P})$,
	i.e. a centered Gaussian process satisfying that $\mathbb{E}[W(h)W(g)]=\langle h,g\rangle_{\mathcal{H}}$.
	We denote by $\mathcal{S}$ the family of smooth random variables,
	that is, the collection of all random variables that can be written as
	\begin{equation}
		F=f(W(h_{1}),\ldots,W(h_{N})),\,\,N\geq1,\label{eq:smooth_rvs_def}
	\end{equation}
	in which $f\in\mathcal{C}_{b}^{\infty}(\mathbb{R})$ (the class of
	real-valued bounded functions with bounded derivatives of any order)
	and $h_{i}\in\mathcal{H}$, for $i=1,\ldots,N$. Given $F\in\mathcal{S}$
	with representation (\ref{eq:smooth_rvs_def}), we define and denote
	the Malliavin derivative of $F$ as 
	\[
	DF:=\sum_{i=1}^{N}\frac{\partial f}{\partial x_{i}}(W(h_{1}),\ldots,W(h_{N}))h_{i}.
	\]
	Note that for every $p\geq1$, $DF$ can be seen as an element of
	$L^{p}(\Omega;\mathcal{H})$ (the class of $\mathcal{H}$-valued random
	mappings $G:\Omega\rightarrow\mathcal{H}$ satisfying that $\mathbb{E}(\rVert G\rVert_{\mathcal{H}}^{p})<\infty$).
	Given $p\geq1$, $\mathbb{D}^{1,p}$ will denote the closure of $\mathcal{S}$
	respect to the norm
	\[
	\rVert F\rVert_{1,p}=\left[\mathbb{E}(\rvert F\rvert^{p})+\mathbb{E}(\rVert DF\rVert_{\mathcal{H}}^{p})\right]^{1/p}.
	\]
	The Malliavin derivative is a closable operator from $L^{p}(\Omega)=L^{p}(\Omega,\mathcal{F},\mathbb{P})$
	into $L^{p}(\Omega;\mathcal{H})$ and the domain of such extension
	is exactly $\mathbb{D}^{1,p}$. The symbol $\delta$ will indicate
	the adjoint operator of $D$ in $L^{2}(\Omega)$, in other words $\delta$ is
	an operator from $\mathrm{Dom}\delta\subseteq L^{2}(\Omega;\mathcal{H})$
	to $L^{2}(\Omega)$ satisfying that 
	\[
	\mathbb{E}(\langle DF,G\rangle_{\mathcal{H}})=\mathbb{E}(F\delta(G)),\,\,\,\forall\,G\in\mathrm{Dom}\delta.
	\]
	Now, let $\mathcal{U}$ be another separable real Hilbert space with
	inner product $\langle,\rangle_{\mathcal{U}}$ and denote by $\mathcal{S}_{\mathcal{U}}$
	the collection of $\mathcal{U}$-valued random mappings that can
	be written as $G=\sum_{i=1}^{N}F_{i}u_{i}$, where $F_{i}\in\mathcal{S}$
	and $u_{i}\in\mathcal{U}$ for $i=1,\ldots,N$ for some $N\in\mathbb{N}$.
	By $\mathbb{D}^{1,p}(\mathcal{U})$ we mean the closure of $\mathcal{S}_{\mathcal{U}}$
	under the norm 
	\[
	\rVert G\rVert_{1,p,\mathcal{U}}:=\left[\mathbb{E}(\rVert G\rVert_{\mathcal{U}}^{p})+\mathbb{E}(\rVert DG\rVert_{\mathcal{H}\otimes\mathcal{U}}^{p})\right]^{1/p};\,\,\,G=\sum_{i=1}^{N}F_{i}u_{i},
	\]
	where $DG:=\sum_{i=1}^{N}DF_{i}\otimes u_{i}$.
	
	The following properties are known:
	\begin{enumerate}
		\item For any $F\in\mathbb{D}^{1,2}$ and $G\in\mathrm{Dom}\delta,$ such
		that $FG\in L^{2}(\Omega;\mathcal{H})$, it holds that 
		\begin{equation}
			F\delta(G)=\delta(FG)+\langle DF,G\rangle_{\mathcal{H}}.\label{eq:product_rule_malliavin}
		\end{equation}
		\item The operator $\delta$ is continuous from $\mathbb{D}^{1,p}(\mathcal{H})$
		into $L^{p}(\Omega)$, put differently, if $G\in\mathbb{D}^{1,p}(\mathcal{H})$
		then $G\in\mathrm{Dom}\delta$ and there is a constant $C_{p}$
		such that
		\begin{equation}
			\mathbb{E}(\rvert\delta(G)\rvert^{p})^{1/p}\leq C_{p}\rVert G\rVert_{1,p,\mathcal{H}},\,\,\,\forall\,G\in\mathbb{D}^{1,p}(\mathcal{H}).\label{eq:meyers_ine}
		\end{equation}
	\end{enumerate}
	We refer the reader to the monograph \cite{Nualart06} for a concise and detailed
	exposition on the theory of Malliavin calculus.
	
	\section{The Euler-Maruyama method\label{sec:The-Euler-Maruyama-method}}
	
	In this part we introduce the Euler-Maruyama method associated to
	the SDDE 
	\begin{equation}
		X_{t}=\begin{cases}
			x_{0}(0)+\int_{0}^{t}\int_{[0,\tau]}b(X_{s-r})\eta(\mathrm{d}r)\mathrm{d}s+B_{t}, & t\geq0,\\
			x_{0}(t), & -\tau\leq t<0,
		\end{cases}\label{SDDE}
	\end{equation}
	where $\eta$ is a finite signed measure on $[0,\tau]$, $\tau>0$,
	$B$ is an fBm, and $x_{0}$ is a continuous and deterministic function
	with the convention that $x_{0}$ vanishes outside of $[-\tau,0]$,
	and investigate the asymptotic behaviour of its error. Existence and
	uniqueness of solutions of (\ref{SDDE}) have been extensively studied
	in the case when $H>1/2$ in the case when $b$ is Lipschitz continuous.
	See for instance \cite{FerrantRovira10,BOUFHajji11,BesaluRov12,CaraballoGarridoTanigu11,LeonTindel12}
	and references therein. Since the noise in (\ref{SDDE}) is additive,
	the arguments used in these works can be extended to the case $H\leq1/2$.
	
	\begin{remark}\label{remarkboundsmoments}The following remarks are
		pertinent: 
		\begin{enumerate}
			\item Any process satisfying (\ref{SDDE}) inherit the path properties of
			$B$. Thus, any solution of such equation has local $\lambda$-Hölder
			continuous paths on $[0,+\infty)$ for every $\lambda<H$.
			\item If $b$ is of linear growth, and $X$ solves (\ref{SDDE}), a localization
			argument along with Gronwall's inequality and the self-similarity
			of $B$ (see also Theorem 1.1 in \cite{NovikovValkeila99}) show
			that for all $p\geq 1$ and every $T>0$ 
			\begin{equation}
				\mathbb{E}(\sup_{0\leq t\leq T}\lvert X_{t}\lvert^{p})\leq C_{1}\mathbb{E}(\sup_{0\leq t\leq T}\lvert B_{t}\lvert^{p})e^{C_{2}T^{p-1}}=C_{1}e^{C_{2}T^{p-1}}T^{pH},\label{eq:bounmomentssol}
			\end{equation}
			for some positive constants $C_{1},C_{2}$ depending only on $T,p,H,\eta,x_{0}$
			and $b$.
			\item Relation (\ref{eq:bounmomentssol}) together with (\ref{eq:HolderexplicitfBm})
			guarantee that 
			\begin{equation}
				\lvert X_{t}-X_{s}\lvert\leq\xi_{\lambda,T}\lvert t-s\lvert^{\lambda},\,\,\forall\,t,s\in[0,T],T>0,\label{eq:Holdercond_X}
			\end{equation}
			where $0<\lambda<H$ and $\xi_{\lambda,T}$ is a positive random variable
			with finite moments of all orders. 
		\end{enumerate}
	\end{remark} 
	
	In what is left of this work, we will always assume that the drift
	component \textbf{$b$} is at least of class $\mathcal{C}^{1}$ and of linear
	growth. Within this framework, the unique solution of (\ref{SDDE})
	will be denoted by $X$. Now, by setting 
	\[
	\mathcal{T}(s):=[s/\Delta_{n}]\Delta_{n},\,s\geq0,\,\,\Delta_{n}:=\tau/n,
	\]
	we define and denote the Euler-Maruyama scheme associated to (\ref{SDDE})
	as the stochastic process 
	\begin{equation}
		X_{t}^{n}:=\begin{cases}
			x_{0}(0)+\int_{0}^{t}\int_{[0,\tau]}b(X_{\mathcal{T}(s)-\mathcal{T}(r)}^{n})\eta(\mathrm{d}r)\mathrm{d}s+B_{t}, & t\geq0,\\
			x_{0}(t), & -\tau\leq t<0.
		\end{cases}\label{eq:Eulermarscheme}
	\end{equation}
	Set $t_{i}=i\Delta_{n}$ for $i=0,1,2,\ldots$, and observe that for
	every $t\in[t_{i-1},t_{i})$ 
	\[
	X_{t}^{n}=X_{t_{i-1}}^{n}+(t-t_{i-1})\left\{ b(X_{t_{i-1}-\tau}^{n})\eta(\{\tau\})+\sum_{j=1}^{n}b(X_{t_{i-1}-t_{j-1}}^{n})\eta([t_{j-1},t_{j}))\right\} +B_{t}-B_{t_{i-1}}.
	\]
	This in particular implies that $X^{n}$ is jointly measurable and
	$\lambda$-Hölder continuous for any $0<\lambda<H$. We are interested
	on the behaviour of the error process 
	\begin{equation}
		U_{t}^{n}=X_{t}-X_{t}^{n},\,\,\,t\geq-\tau.\label{eq:deferror}
	\end{equation}
	A simple application of Gronwall's inequality, relation \eqref{eq:bounmomentssol},
	and the self-similarity of $B$ result in
	\begin{equation}
		\mathbb{E}(\sup_{t\leq T}\vert U_{t}^{n}\vert^{p})\leq C\Delta_{n}^{pH},\,\,p\geq1,\label{eq:naiveineq}
	\end{equation}
	for some contant $C$ independent of $n$. 
	
	\begin{remark}\label{rmk_momentEuler}Note that if $b$ is of linear
		growth, then 
		\[
		\rvert X_{t}\rvert\leq C(\sup_{s\leq T}\vert U_{s}^{n}\vert+\sup_{s\leq T}\vert X_{s}\vert)+\rvert B_{t}\rvert\,\forall\,t\in[0,T],T>0.
		\]
		In particular, for any $p\geq1$
		\begin{equation}
			\mathbb{E}(\sup_{0\leq t\leq T}\lvert X_{t}^{n}\lvert^{p})<\infty,\label{eq:bounmoment_Euler}
		\end{equation}
		thanks to (\ref{eq:naiveineq}) and (\ref{eq:bounmomentssol}).
		
	\end{remark} The rate $\Delta_{n}^{H}$ is by no means optimal. In
	fact, \cite{HuLiuNualart16} have shown that when $\eta(\mathrm{d}r)$
	is the Dirac's delta measure at $0$ (i.e. no delay) and $H>3/4$,
	the optimal rate is $\Delta_{n}$. We will see that if $\eta(\mathrm{d}r)$
	is continuous, this is not in general the case.
	
	Before presenting our main findings, let us introduce some notation.
	Given a measurable process $(N_{t})_{t\geq0}$, $(U(N)_{t})_{t\geq0}$
	will denote the unique solution (if it exists) of the semilinear SDDE
	\begin{equation}
		U(N)_{t}=\begin{cases}
			\int_{0}^{t}\int_{[0,s]}b^{\prime}(X_{s-r})U_{s-r}\eta(\mathrm{d}r)\mathrm{d}s+N_{t}, & t\geq0,\\
			0, & -\tau\leq t<0.
		\end{cases}\label{eq:semilinear limit}
	\end{equation}
	Furthermore, if $b^{\prime}$ is bounded, $R$ will represent the
	differential resolvent (see Subsection \ref{subsec:Differential-resolvents})
	of 
	\[
	\kappa(t,A)=-\int_{[0,t]}b^{\prime}(X_{t-r})\mathbf{1}_{A}(t-r)\eta(\mathrm{d}r),\,\,t\geq0,\,\,A\in\mathcal{B}([0,+\infty)).
	\]
	We remind to the reader that we are assuming that $b$ is of linear
	growth. 
	
	\begin{theorem}[The non-rough case]\label{thmerrordisGaussiannonrough}
		Suppose that $b$ is of class $\mathcal{C}^{2}$ with $b^{\prime}\in\mathcal{C}_{b}^{1}$
		and that $x_{0}$ is of class $\mathcal{C}^{1}$. If $H\geq1/2$,
		then for every $t>0$ and $p\geq1$ there is $C>0$ such that 
		\begin{equation}
			\sup_{0\leq t\leq T}\mathbb{E}(\vert U_{t}^{n}\vert^{p})^{1/p}\leq C\Delta_{n}.\label{eq:optimratenonrough}
		\end{equation}
		Furthermore:
		\begin{enumerate}
			\item If $H>1/2$, it holds that 
			\begin{equation}
				\frac{1}{\Delta_{n}}(U^{n}-\mathfrak{B}^{n})\overset{u.c.p}{\longrightarrow}U(N),\label{eq:Lpconv}
			\end{equation}
			where 
			\[
			N_{t}=\frac{1}{2}\int_{[0,\tau]}(b(X_{t-r})-b(X_{-r}))\eta(\mathrm{d}r),\,\,t\geq0,
			\]
			and 
			\begin{equation}
				\mathfrak{B}_{t}^{n}=\int_{0}^{t}\int_{[0,\tau]}R(t,s)[b(X_{s-r})-b(X_{s-\mathcal{T}(r)})]\eta(\mathrm{d}r)\mathrm{d}s,\,\,t\geq0.\label{eq:Bias}
			\end{equation}
			\item If $H=1/2$, as $n\rightarrow\infty$
			\begin{equation}
				\frac{1}{\Delta_{n}}(U^{n}-\mathfrak{B}^{n})\overset{\mathcal{F}\text{-}C[0,T]}{\Longrightarrow}U(N),\label{eq:Lpconv-1}
			\end{equation}
			in which 
			\[
			N_{t}=\frac{1}{2}\int_{[0,\tau]}(b(X_{t-r})-b(X_{-r}))\eta(\mathrm{d}r)+\frac{1}{\sqrt{12}}\int_{[0,t]}\int_{0}^{t-r}b^{\prime}(X_{s})\mathrm{d}\tilde{W}_{s},\,\,t\geq0,
			\]
			with $\tilde{W}$ a Brownian motion defined on an extension of $(\Omega,\mathcal{F},\mathbb{P})$
			which is in turn independent of $\mathcal{F}$.
		\end{enumerate}
	\end{theorem}
	
	\begin{remark}\label{nonconvfV}It is not difficult to see from our
		conditions on $b$ that for all $p\geq1$ and $T>0$
		\[
		\sup_{0\leq t\leq T}\vert\mathfrak{B}_{t}^{n}\vert\leq\zeta_{T}\Delta_{n},
		\]
		and $\zeta_{T}$ is a positivse random variable with finite moments
		of all orders. However, in general $\frac{1}{\Delta_{n}}\mathfrak{B}^{n}$
		does not fulfill a limit theorem. For instance, if $\eta$ is purely
		atomic and $r$ is an atom such that $\mathcal{T}(r)\neq r$. Despite
		of this negative result, if $\eta$ admits a continuous density then
		$\frac{1}{\Delta_{n}}\mathfrak{B}^{n}$ converges. Specifically, for
		any $H\in(0,1)$
		\begin{equation}
			\frac{1}{\Delta_{n}}\mathfrak{B}^{n}\overset{u.c.p}{\longrightarrow}-U(N^{\prime}),\label{eq:limitforVn}
		\end{equation}
		where $N_{t}^{\prime}=\frac{1}{2}\int_{[0,\tau]}(b(X_{t-r})-b(X_{-r}))\eta(\mathrm{d}r).$
		Thus, 
		\begin{enumerate}
			\item For $H>1/2$, $\frac{1}{\Delta_{n}}U^{n}\overset{u.c.p}{\longrightarrow}0$. 
			\item If $H=1/2$, $\frac{1}{\Delta_{n}}U^{n}\overset{\mathcal{F}\text{-}C[0,T]}{\Longrightarrow}U(\tilde{N}),$where
			\[
			\tilde{N}_{t}=\frac{1}{\sqrt{12}}\int_{[0,t]}\int_{0}^{t-r}b^{\prime}(X_{s})\mathrm{d}\tilde{W}_{s},\,\,t\geq0.
			\]
		\end{enumerate}
		A proof for (\ref{eq:limitforVn}) will be presented in Section \ref{subsec:Proof-of-Theoremnonrough}.
		
	\end{remark}
	
	\begin{remark}\label{rmk_result_nodelay}Let us see how our results
		compare with those obtained in \cite{HuLiuNualart16}. If $H>1/2$
		and $\eta(\mathrm{d}r)$ is the Dirac's delta measure at $0$, then
		$V^{n}\equiv0$. Furthermore, path-wise integration with respec to
		$X$ is possible. Hence, by the change of variables formula for the
		Riemann-Stieltjes integral (see for instance \cite{Zale98}, Theorem
		4.3.1) and (\ref{SDDE}) we also get that
		\begin{align*}
			N_{t}= & \frac{1}{2}\int_{0}^{t}b^{\prime}(X_{s})\mathrm{d}X_{s},\,\,\,t\geq0.\\
			= & \frac{1}{2}\int_{0}^{t}b^{\prime}(X_{s})b(X_{s})\mathrm{d}s+\frac{1}{2}\int_{0}^{t}b^{\prime}(X_{s})\mathrm{d}B_{s}.
		\end{align*}
		In consequence, the process $U=U(N)$ in Theorem \ref{nonconvfV}
		satisfy the semilinear stochastic differential equation 
		\[
		U_{t}=\int_{0}^{t}b^{\prime}(X_{s})U_{s}\mathrm{d}s+\frac{1}{2}\int_{0}^{t}b^{\prime}(X_{s})b(X_{s})\mathrm{d}s+\frac{1}{2}\int_{0}^{t}b^{\prime}(X_{s})\mathrm{d}B_{s},\,\,\,t\geq0.
		\]
		which is exactly the conclusion of Theorem 8.1 in \cite{HuLiuNualart16}.
		
	\end{remark}
	
	Our discussion in Remark \ref{nonconvfV} demonstrates that when $H>1/2$
	and the delay $\eta$ is absolutely continuous, the rate $\Delta_{n}$ --
	in constrast with the non-delayed case -- is suboptimal. Motivated
	by this, in the next result we further investigate how the rate can
	be improved in such a situation.
	
	\begin{theorem}[The non-rough case]\label{thmerrordisGaussiannonrough-1}
		Let the assumptions of Theorem \ref{thmerrordisGaussiannonrough}
		hold and let $H>1/2$. Suppose that $\eta$ admits a Lipschitz continuous
		density and that $x_{0}$ is of class $\mathcal{C}^{2}$. Then, for
		every $t>0$ and $p\geq1$ there is $C>0$ such that 
		\begin{equation}
			\sup_{0\leq t\leq T}\mathbb{E}(\vert U_{t}^{n}\vert^{p})^{1/p}\leq C\Delta_{n}^{H+1/2}.\label{eq:optimratenonrough-1}
		\end{equation}
		
	\end{theorem}
	
	Not surprisingly, the rough case is subtantially more delicated. In
	this situation we managed to identify only the rate of convergence.
	
	\begin{theorem}[The rough case]\label{themroughcase} Suppose that
		$b$ is of class $\mathcal{C}^{3}$ with $b^{\prime}\in\mathcal{C}_{b}^{2}$
		and that $x_{0}$ is of class $\mathcal{C}^{1}$. If $0<H<1/2$, then
		for every $T>0$ and $p\geq1$ there is $C>0$ such that 
		\begin{equation}
			\sup_{0\leq t\leq T}\mathbb{E}(\vert U_{t}^{n}\vert^{p})^{1/p}\leq C\Delta_{n}^{3H\land(H+1/2)}.\label{eq:optimraterough}
		\end{equation}
		
	\end{theorem}
	
	\begin{remark}Our proof suggests (see Remarks \ref{rmk_bettererate_rough1})
		the following:
		\begin{enumerate}
			\item When $b$ is of class $\mathcal{C}^{4}$ with $b^{\prime}\in\mathcal{C}_{b}^{3}$
			the rate $\Delta_{n}^{4H\land(H+1/2)}$ can be attained. Based on
			this, we conjecture that for any $H\in(0,1)$ the optimal rate is
			$\Delta_{n}^{H+1/2}$ whenever $b$ is of class $\mathcal{C}^{\infty}$
			with $b^{\prime}\in\mathcal{C}_{b}^{\infty}$. In an earlier version
			of this work (see also our companion paper \cite{Sauri25}), we prove that when $b$ is linear this is indeed
			the case.
			\item \eqref{eq:optimraterough} can be slightly improved. Indeed, under the assumptions of Theorem
			\ref{themroughcase}, for every $T>0$ and $p\geq1$, we can find
			a constant that does not depend on $n$ such that
			\[
			\mathbb{E}(\vert U_{t}^{n}-U_{v}^{n}\vert^{p})^{1/p}\leq C\Delta_{n}^{3H\land(H+1/2)}(t-v)^{1/2},\,\, 0\leq v\leq t \leq T. 
			\]
			This tightness result hints that, as in the case $H=1/2$, $U^{n}$
			converges (after a normalization) towards an SDDE driven by a functional of a standard
			Brownian. This is proved in our companion paper \cite{Sauri25}
			when $b$ is linear.
		\end{enumerate}
	\end{remark}
	
	\section{Proofs\label{sec:Proofs}}
	
	Throughout all our proofs, the non-random positive constants (independent of $n$) will
	be denoted by the generic symbol $C>0$, and they may change from
	line to line. As it is costumary, we will write $x\lesssim y$ whenever $x\leq C y$ where the constant $C$ being independent of $n$. 
	We recall to the reader that we are always assuming
	that $b$, the drift component of our SDDE, is of class $\mathcal{C}^{1}$ and
	of linear growth.
	
	We will often use the following notation: Given a sequence of processes
	$(H_{t}^{n})_{0\leq t\leq T}$ and a deterministic sequence $r_{n}$,
	we will write $H_{t}^{n}=\mathscr{O}_{p}(r_{n})$ if there is a constant
	$C>0$ independent of $n$ such that
	\[
	\sup_{0\leq t\leq T}\mathbb{E}(\rvert H_{t}^{n}\rvert^{p})^{1/p}\leq Cr_{n},
	\]
	Similarly, we use the notation $H_{t}^{n}=\mathscr{O}_{p}^{u}(r_{n})$ whenever
	\[
	\mathbb{E}(\sup_{0\leq t\leq T}\rvert H_{t}^{n}\rvert^{p})^{1/p}\leq Cr_{n}.
	\]

	\subsection{Measurability of the differential resolvent}
	
	Recall that a Volterra kernel measure of type $B^{\infty}$ on $J\subseteq[0,+\infty)$,
	is a mapping $\mu:J\times\mathcal{B}(J)\rightarrow\mathbb{R}$ such
	that: $i)$ $\mu(t,\cdot)$ is signed finite measure with support
	on $[0,t]$. $ii)$ For every Borelian set $A$, $\mu(\cdot,A)$ is
	measurable. $iii)$ The mapping $t\in J\mapsto\lvert\mu\lvert(t,J)$
	is bounded. The collection of all Volterra kernel measure of type
	$B^{\infty}$ on $J$ will be denoted as $\mathcal{M}(B^{\infty},J)$.
	Note that if $\mu\in\mathcal{M}(B^{\infty},J)$, then $\mu\in\mathcal{M}(B^{\infty},I)$,
	for every $I\subseteq J$. Reciprocally, if $\mu\in\mathcal{M}(B^{\infty},I)$
	and $I\subseteq J$, we can lift $\mu$ as an element of $\mathcal{M}(B^{\infty},J)$
	by letting 
	\begin{equation}
		\mu(t,A)=\begin{cases}
			\mu(t,A\cap I) & ,\,\,t\in I,A\in\mathcal{B}(J);\\
			0 & \text{otherwise}.
		\end{cases}\label{eq:liftmu_in_I_to_J}
	\end{equation}
	Let us now introduce some notation. Given $\mu,\nu\in\mathcal{M}(B^{\infty},J)$,
	we set 
	\[
	\mu\star\nu(t,A):=\int_{J}\nu(s,A)\mu(t,\mathrm{d}s),\,\,t\in J,A\in\mathcal{B}(J),
	\]
	and if $k:[0,+\infty)^{2}\rightarrow\mathbb{R}$ is a Volterra kernel
	($k(t,s)=0$ if $s>t$) measurable and bounded, we define 
	\[
	k\star\mu(t,A):=\int_{J}k(t,s)\mu(s,A)\mathrm{d}s,\,\,t\in J,A\in\mathcal{B}(J),
	\]
	and 
	\[
	\mu\star k(t,s):=\int_{J}k(u,s)\mu(t,\mathrm{d}u),\,\,s,t\in J.
	\]
	It is not difficult to see that $\mu\star\nu$, $k\star\mu\in(B^{\infty},J)$
	and $\mu\star k$ are measurable and bounded Volterra kernels. The
	space $\mathcal{M}(B^{\infty},J)$ is a Banach algebra with multiplication
	$\star$ (see \cite{GripenbergLondenStaffans90}, Theorem 2.3) if
	endowed with the norm 
	\[
	\parallel\mu\parallel_{\mathcal{M}_{\infty}(J)}=\sup_{t\in J}\lvert\mu\lvert(t,J).
	\]
	Lastly, if $\mu,\rho\in\mathcal{M}(B^{\infty},J)$, we say that $\rho$
	is a \textit{measure resolvent} of $\mu$, if 
	\[
	\rho+\mu\star\rho=\rho+\rho\star\mu=\mu.
	\]
	Let us now briefly summarize the construction given in \cite{GripenbergLondenStaffans90},
	Chapter 10, of the differential resolvent of $\mu\in\mathcal{M}(B^{\infty},[0,\infty))$.
	Put $k(t,s)=\mathbf{1}_{s\leq t}$ and set 
	\[
	\nu(t,A)=k\star\mu(t,A)=\int_{[0,t]}\mu(s,A)\mathrm{d}s.
	\]
	If $\nu$ admits a measure resolvent, say $\rho$, then the unique
	differential resolvent of $\mu$ is given by 
	\begin{equation}
		r(t,s)=k(t,s)-\rho\star k(t,s)=\mathbf{1}_{s\leq t}-\rho(t,[s,t]).\label{eq:diffresolvconnst}
	\end{equation}
	We now give a local representation of $\rho$. Let $[0,T^{\prime})$
	be such that $s,t\in[0,T^{\prime})$ and $0=T_{0}<T_{1}<\cdots T_{N}=T^{\prime}$
	be a finite subdivision of $[0,T^{\prime})$, such that $\parallel\nu\parallel_{\mathcal{M}_{\infty}([T_{i-1},T_{i}))}<1$,
	for all $i=1,\ldots,N.$ By the Banach fixed-point Theorem, there
	is a unique $\bar{\rho}_{i}\in\mathcal{M}(B^{\infty},[T_{i-1},T_{i}))$
	resolvent for $\nu$ as element of $\mathcal{M}(B^{\infty},[T_{i-1},T_{i}))$
	(i.e. $\nu$ restricted to $[T_{i-1},T_{i})$) which can be written
	as 
	\begin{equation}
		\bar{\rho}_{i}=\lim_{n\rightarrow\infty}\bar{\rho}_{i}^{(n)}=\lim_{n\rightarrow\infty}\sum_{j=1}^{n+1}(-1)^{j-1}\nu^{\star j},\,\,\,\nu^{\star j}:=\underbrace{\nu\star\nu\star\cdots\star\nu}_{j\text{-times}},\label{eq:rho_hatdef}
	\end{equation}
	where the limit is taken w.r.t. to $\parallel\cdot\parallel_{\mathcal{M}_{\infty}([T_{i-1},T_{i}))}$.
	$\rho$ is now constructed iteratively. First, for $i=1,\ldots,n$,
	we build $\rho_{i}\in\mathcal{M}(B^{\infty},[0,T_{i}))$ (based on
	$\bar{\rho}_{i}$) as the measure resolvent of $\nu$ as element of
	$\mathcal{M}(B^{\infty},[0,T_{i}))$. For $i=1$, we set $\rho_{1}=\bar{\rho}_{1}$.
	For $i=2$, first we lift $\rho_{1}$ and $\bar{\rho}_{2}$ as elements
	of $\mathcal{M}(B^{\infty},[0,T_{2}))$ via the procedure given in
	(\ref{eq:liftmu_in_I_to_J}). Then, we define 
	\[
	\mathcal{M}(B^{\infty},[0,T_{2}))\ni\rho_{2}=\nu-\nu\star\rho_{1}-\bar{\rho}_{2}\star\nu+\bar{\rho}_{2}\star\nu\star\rho_{1}.
	\]
	Note that 
	\[
	\rho_{2}=\lim_{n\rightarrow\infty}\rho_{2}^{(n)}:=\lim_{n\rightarrow\infty}(\nu-\nu\star\rho_{1}^{(n)}-\bar{\rho}_{2}^{(n)}\star\nu+\bar{\rho}_{2}^{(n)}\star\nu\star\rho_{1}^{(n)}),\,\,\,\rho_{1}^{(n)}=\bar{\rho}_{1}^{(n)},
	\]
	on $\mathcal{M}(B^{\infty},[0,T_{2}))$. For $i=3,\ldots,N$ we repeat
	this procedure, i.e. we first extend $\rho_{i-1}$ and $\bar{\rho}_{i}$
	as elements of $\mathcal{M}(B^{\infty},[0,T_{i}))$ and then set $\rho_{i}\in\mathcal{M}(B^{\infty},[0,T_{i}))$
	as 
	\[
	\rho_{i}=\nu-\nu\star\rho_{i-1}-\bar{\rho}_{i}\star\nu+\bar{\rho}_{i}\star\nu\star\rho_{i-1}.
	\]
	We also have that 
	\begin{equation}
		\rho_{i}=\lim_{n\rightarrow\infty}\rho_{i}^{(n)}:=\lim_{n\rightarrow\infty}(\nu-\nu\star\rho_{i-1}^{(n)}-\bar{\rho}_{i}^{(n)}\star\nu+\bar{\rho}_{i}^{(n)}\star\nu\star\rho_{i-1}^{(n)}),\label{eq:rho_i_seq}
	\end{equation}
	on $\mathcal{M}(B^{\infty},[0,T_{i}))$. Finally, we set $\rho=\rho_{N}$
	and note that $\rho$ coincides with $\rho_{i}$ on $\mathcal{M}(B^{\infty},[0,T_{i}))$.
	
	\begin{proof}[Proof of Proposition \ref{propmeasurableresolv}]
		
		During the proof we will keep the notation introduced above. Fix $T>0$
		and choose $T^{\prime}>T$, such that (\ref{eq:rho_hatdef})-(\ref{eq:rho_i_seq})
		hold. Our proof relies on the following fact (which is a simple consequence
		of (\ref{eq:diffresolvconnst}), that $\rho$ coincides with $\rho_{i}$
		on $\mathcal{M}(B^{\infty},[0,T_{i}))$, and that $\rho_{i}\rightarrow\rho_{i}^{(n)}$
		on $\mathcal{M}_{\infty}([0,T_{i}))$): If $r$ is the differential
		resolvent of $\mu\in\mathcal{M}(B^{\infty},[0,\infty))$, then 
		\begin{equation}
			r(t,s)=\lim_{n\rightarrow\infty}r_{n}(t,s),\,\,\,s,t\in[0,T_{i}),\label{eq:limitdifresolv}
		\end{equation}
		where 
		\[
		r_{n}(t,s)=k(t,s)-\rho_{i}^{(n)}\star k(t,s)=k(t,s)-\rho_{i}^{(n)}(t,[s,t]),
		\]
		where $\rho_{i}^{(n)}$ as in (\ref{eq:rho_i_seq}). In view of this,
		it is enough to show that the mapping $(t,s,\omega)\in[0,T_{i})^{2}\times\Omega\mapsto\rho_{i}^{(n)}(t,[s,t],\omega)$
		is measurable, for every $i=1,\ldots,N$. Here $\rho_{i}^{(n)}(\cdot,\omega)$
		is as in (\ref{eq:rho_i_seq}) with $\mu=\kappa(\cdot,\omega)$. Put
		$\mathscr{U}_{i}:=\{(t,s)\in[0,T_{i})^{2}:s\leq t\}$. Since for $s>t$,
		$\rho_{i}^{(n)}(t,[s,t])$ does not depend on $\omega$, we only need
		to concentrate on the case in which $(s,t)\in\mathscr{U}_{i}$. To
		see that this is the case we first show that for all $j=1,\ldots$
		the mapping $(t,(b,a),\omega)\in[0,T_{i})\times\mathscr{U}_{i}\times\Omega\mapsto\nu^{\star j}(t,[a,b],\omega)$
		is measurable. We proceed by induction. The validity of our claim
		for $j=1$ follows from the fact that
		\[
		\nu(t,[a,b],\omega)=-\int_{0}^{\infty}\int_{0}^{\infty}\mathbf{1}_{[0,t]}(u)\mathbf{1}_{[0,+\infty)}(u-r)b_{u-r}(\omega)\mathbf{1}_{[a,b]}(u-r)\eta(\mathrm{d}r)\mathrm{d}u,\,\,(t,(b,a),\omega)\in[0,T_{i})\times\mathscr{U}_{i}\times\Omega,
		\]
		and Fubini's Theorem. Now, suppose that $(t,(b,a),\omega)\in[0,T_{i})\times\mathscr{U}_{i}\times\Omega\mapsto\nu^{\star j}(t,[a,b],\omega)$
		is measurable for $j=k$. In view that 
		\[
		\nu^{\star k+1}(x,[a,b],\omega)=-\int_{0}^{\infty}\int_{0}^{\infty}\mathbf{1}_{[0,t]}(u)\mathbf{1}_{[0,+\infty)}(u-r)b_{u-r}(\omega)\nu^{\star j}(u-r,[a,b],\omega)\eta(\mathrm{d}r)\mathrm{d}u,
		\]
		we have once again, due to Fubini's Theorem, that our claim holds.
		We deduce from this property and (\ref{eq:rho_hatdef}) that the function
		$(t,s,\omega)\in\mathscr{U}_{i}\times\Omega\mapsto\bar{\rho}_{i}^{(n)}(t,[s,t],\omega)$
		is measurable. Applying this into (\ref{eq:rho_i_seq}) along with
		another induction argument let us conclude that $(t,s,\omega)\in[0,T_{i})^{2}\times\Omega\mapsto\rho_{i}^{(n)}(t,[s,t],\omega)$
		is indeed measurable, as required.\end{proof}
	
	\subsection{Malliavin differentiability of the solution of (\ref{SDDE})}
	
	Multiple estimates used in our proof of Theorems (\ref{thmerrordisGaussiannonrough})-(\ref{themroughcase})
	relay in (\ref{eq:meyers_ine}). Therefore, in this part we verify
	that within our framework the solution to the SDDE (\ref{SDDE}) is
	always differentiable in the Malliavin sense. Fix $T>0$. Given the
	representation of $B$ in terms of the two-sided Wiener process $W$
	(see \eqref{fBmdef}), we will consider the Hilbert space $\mathcal{H}=L^{2}((-\infty,T],\mathrm{d}s;\mathbb{R})$,
	with inner product $\langle h,g\rangle_{\mathcal{H}}=\int_{-\infty}^{T}h(s)g(s)\mathrm{d}s$,
	so that our basic isonormal Gaussian process becomes $\{\int_{-\infty}^{T}h(s)\mathrm{d}W_{s}:h\in\mathcal{H}\}$.
	Without loss of generality we may and do assume that $\mathcal{F}=\mathcal{F}_{T}^{W}$,
	where $\mathcal{F}_{T}^{W}$ is the completion of the $\sigma$-algebra
	generated by $(W_{t}:t\leq T)$. 
	
	\begin{proposition} \label{prop_Malliavinder}If $b^{\prime}\in\mathcal{C}_{b}$,
		then $X_{t}\in\mathbb{D}^{1,p}$ for all $t\in[-\tau,T]$ and any
		$p>1$. Furthermore, there is a version of $(u,t,\omega)\in(-\infty,T]\times[-\tau,T]\times\Omega\mapsto D_{u}X_{t}(\omega)$
		that is jointly measurable and satisfies the semi-linear SDDE
		\begin{equation}
			D_{u}X_{t}=\begin{cases}
				\int_{0}^{t}\int_{[0,s]}b^{\prime}(X_{s-r})D_{u}X_{s-r}\eta(\mathrm{d}r)\mathrm{d}s+K(t,u), & t>0,u\leq t;\\
				0 & -\tau\leq t\leq0;\text{or if }u>t;
			\end{cases}\label{Malliav_Der_SDDE}
		\end{equation}
		and $D_{u}X_{t}=0$ for $u>t$, where
		\begin{equation}
			K(t,u)=(t-u)_{+}^{H-1/2}-(-u)_{+}^{H-1/2}.\label{eq:kernel_fBm}
		\end{equation}
	\end{proposition}
	
	\begin{proof}Standard arguments (see for instance \cite{KaratzasShreve91},
		Theorem 5.2.9) show that the Piccard sequence
		\[
		X_{t}^{(n)}=\begin{cases}
			x_{0}(0)+\int_{0}^{t}\int_{[0,\tau]}b(X_{s-r}^{(n-1)})\eta(\mathrm{d}r)\mathrm{d}s+Z_{t}, & t\geq0,\\
			x_{0}(t) & -\tau\leq t<0,
		\end{cases}
		\]
		where $X_{t}^{(0)}=x_{0}(t\land0)$, for $t\geq-\tau$, and $n\in\mathbb{N}$,
		converges to $X$ in $L^{p}(\mathbb{P})$, $p\geq1$, and uniformly
		in compacts almost surely. Next, to verify that (\ref{Malliav_Der_SDDE})
		holds, we argue as in \cite{FerrantRovira06} (c.f. Theorem 2.2.1
		in \cite{Nualart06}). We are going to show first that, for all $n\in\mathbb{N}$,
		\begin{enumerate}
			\item For every $p>1$, $X_{t}^{(n)}\in\mathbb{D}^{1,p}$, $t\in[-\tau,T].$
			\item There is a version of $(u,t,\omega)\in(-\infty,T]\times[-\tau,T]\times\Omega\mapsto D_{u}X_{t}^{(n)}(\omega)$
			that is measurable.
			\item There is a constant $C_{n}>0$, only depending on $T$ and possibly
			on $n$, such that
			\[
			\sup_{t\in[-\tau,T]}\mathbb{E}(\rVert DX_{t}^{(n)}\rVert_{\mathcal{H}}^{p})\leq C_{n}.
			\]
		\end{enumerate}
		Note that from Proposition 1.5.5 in \cite{Nualart06} and the fact
		that $\mathbb{D}^{1,2}\subseteq\mathbb{D}^{1,p}$ for all $p\in(1,2)$,
		1. is obtained as long as $X_{t}^{(n)}\in\mathbb{D}^{1,2}$ and for
		$p>2$
		\[
		\mathbb{E}(\rVert DX_{t}^{(n)}\rVert_{\mathcal{H}}^{p})=\mathbb{E}\left[\left(\int_{-\infty}^{T}\rvert D_{u}X_{t}^{(n)}\rvert^{2}\mathrm{d}u\right)^{p/2}\right]<\infty.
		\]
		Let us proceed by induction. For $n=1$, 1. and 2. obviously holds
		if $t<0$, so assume that $t\geq0$. Then,
		\[
		X_{t}^{(1)}=x_{0}(0)+\int_{0}^{t}\int_{[0,\tau]}b(x_{0}((s-r)\land0))\eta(\mathrm{d}r)\mathrm{d}s+Z_{t},
		\]
		thus $X_{t}^{(1)}\in\mathscr{F}_{t}^{B}$, $X_{t}^{(1)}\in\mathbb{D}^{1,2}$
		and 
		\[
		D_{u}X_{t}^{(1)}=K(t,u),
		\]
		reason why 
		\[
		\mathbb{E}\left[\left(\int_{-\infty}^{T}\rvert D_{u}X_{t}^{(1)}\rvert^{2}\mathrm{d}u\right)^{p/2}\right]=t^{pH}\leq T^{pH}.
		\]
		We have therefore shown that 1.-3. holds for $n=1$. Suppose now that
		our desired properties hold for $n=k$. Let us see that they are also
		valid fof $k+1$. By definition, 
		\[
		X_{t}^{(k+1)}=x_{0}(0)+\int_{0}^{t}\int_{[0,\tau]}b(X_{s-r}^{(k)})\eta(\mathrm{d}r)\mathrm{d}s+Z_{t}.
		\]
		Let us now see that $F_{t}:=\int_{0}^{t}\int_{[0,\tau]}b(X_{s-r}^{(k)})\eta(\mathrm{d}r)\mathrm{d}s$
		belongs to $\mathbb{D}^{1,2}$ for $t\in[0,T]$ and 
		\begin{equation}
			D_{u}F_{t}=\int_{0}^{t}\int_{[0,\tau]}b^{\prime}(X_{s-r}^{(k)})D_{u}X_{s-r}^{(k)}\eta(\mathrm{d}r)\mathrm{d}s.\label{eq:Deriv_inside_int}
		\end{equation}
		If this were true, we would have that $X_{t}^{(k+1)}\in\mathbb{D}^{1,2}$
		and
		\begin{equation}
			D_{u}X_{t}^{(k+1)}=\int_{0}^{t}\int_{[0,\tau]}b^{\prime}(X_{s-r}^{(k)})D_{u}X_{s-r}^{(k)}\eta(\mathrm{d}r)\mathrm{d}s+K(t,u).\label{eq:Derivative_Piccard}
		\end{equation}
		Hence, by Fubini's Theorem and our induction hypothesis, we would
		have that property 2. above holds for $n=k+1$. Furthermore, by Jensen's
		inequality, for all $p\geq2$ 
		\begin{equation}
			\begin{aligned}\mathbb{E}(\rVert DX_{t}^{(k+1)}\rVert_{\mathcal{H}}^{p})\leq & 2^{p}(T^{pH}+C_{k}\mathfrak{K}_{p}T\times\rvert\eta\rvert([0,\tau])),\end{aligned}
			\label{eq:Exp_it_der}
		\end{equation}
		where $\mathfrak{K}_{p}=\rVert b^{\prime}\rVert_{\infty}^{p}\left(T\times\rvert\eta\rvert([0,\tau]))\right)^{p-1}$,
		which would show that 1. and 3. are also valid for $k+1$, concluding
		this the induction argument. Let us see that $F_{t}\in\mathbb{D}^{1,2}$.
		By the continuity of $X^{(k)}$ and the fact that for all $s,t\in[0,T]$
		\[
		\mathbb{E}(\rvert X_{t}^{(k)}-X_{s}^{(k)}\rvert^{p})\leq C_{p,T}\rvert t-s\rvert^{pH},\,\,\,p\geq1,
		\]
		we conclude that as $m\rightarrow\infty$
		\[
		F_{t}^{m}:=\int_{0}^{t}\int_{[0,\tau]}b(X_{[sm]/m-[rm]/m}^{(k)})\eta(\mathrm{d}r)\mathrm{d}s\overset{L^{p}}{\rightarrow}F_{t}.
		\]
		By the linearity of $D$ and the induction hypothesis we further have
		that $F_{t}\in\mathbb{D}^{1,2}$ and
		\[
		D_{u}F_{t}^{m}=\int_{0}^{t}\int_{[0,\tau]}b^{\prime}(X_{[sm]/m-[rm]/m}^{(k)})D_{u}(X_{[sm]/m-[rm]/m}^{(k)})\eta(\mathrm{d}r)\mathrm{d}s,
		\]
		so that
		\[
		\mathbb{E}(\rVert DF_{t}^{m}\rVert_{\mathcal{H}}^{2})\leq C\sup_{t\in[-\tau,T]}\mathbb{E}(\rVert DX_{t}^{(k)}\rVert_{\mathcal{H}}^{2})\leq C_{k}.
		\]
		Lemma 1.2.3 in \cite{Nualart06}, implies now that $F_{t}\in\mathbb{D}^{1,2}$
		and that for any process $(G_{t})_{t\leq T}$ with $\mathbb{E}(\rVert G\rVert_{\mathcal{H}}^{2})<\infty$,
		we have that as $m\rightarrow\infty$
		\[
		\mathbb{E}(\langle DF_{t}^{m},G\rangle_{\mathcal{H}})\rightarrow\mathbb{E}(\langle DF_{t},G\rangle_{\mathcal{H}}).
		\]
		Thus, to see that (\ref{eq:Deriv_inside_int}) holds, we only need
		to check that 
		\begin{equation}
			\mathbb{E}(\langle DF_{t}^{m},G\rangle_{\mathcal{H}})\rightarrow\mathbb{E}\left[\int_{-\infty}^{T}\left(\int_{0}^{t}\int_{[0,\tau]}b^{\prime}(X_{s-r}^{(k)})D_{u}X_{s-r}^{(k)}\eta(\mathrm{d}r)\mathrm{d}s\right)G_{u}\mathrm{d}u\right].\label{eq:indetification_Deriv}
		\end{equation}
		The Cauchy-Schwartz inequality and the induction hypothesis allow
		us to apply Fubini's Theorem and obtain that
		\[
		\mathbb{E}(\langle DF_{t}^{m},G\rangle_{\mathcal{H}})=\int_{0}^{t}\int_{[0,\tau]}\mathbb{E}\left[\int_{-\infty}^{T}D_{u}\left(b(X_{[sm]/m-[rm]/m}^{(k)})\right)G_{u}\mathrm{d}u\right]\eta(\mathrm{d}r)\mathrm{d}s.
		\]
		Similar arguments to those applied to $F_{t}$ allow us to conclude
		that as $m\rightarrow\infty$
		\[
		\mathbb{E}\left[\int_{-\infty}^{T}D_{u}\left(b(X_{[sm]/m-[rm]/m}^{(k)})\right)G_{u}\mathrm{d}u\right]\rightarrow\mathbb{E}\left[\int_{-\infty}^{T}D_{u}\left(b(X_{s-r}^{(k)})\right)G_{u}\mathrm{d}u\right].
		\]
		Another application of the Cauchy-Schwartz inequality and the induction
		hypothesis enable us the use of the Dominated Convergence Theorem
		to conclude that 
		\begin{align*}
			\mathbb{E}(\int_{-\infty}^{T}D_{u}F_{t}^{m}G_{u}\mathrm{d}u)\rightarrow & \int_{0}^{t}\int_{[0,\tau]}\mathbb{E}\left[\int_{-\infty}^{T}D_{u}\left(b^{\prime}(X_{s-r}^{(k)})\right)G_{u}\mathrm{d}u\right]\eta(\mathrm{d}r)\mathrm{d}s\\
			& =\mathbb{E}\left[\int_{-\infty}^{T}\left(\int_{0}^{t}\int_{[0,\tau]}b^{\prime}(X_{s-r}^{(k)})D_{u}X_{s-r}^{(k)}\eta(\mathrm{d}r)\mathrm{d}s\right)G_{u}\mathrm{d}u\right],
		\end{align*}
		just as needed. Now, iterating (\ref{eq:Exp_it_der}), we conclude
		that for all $n\in\mathbb{N}$ there is $\mathfrak{K}_{p}$ independent
		of $n$, such that 
		\begin{align*}
			\mathbb{E}(\rVert DX_{t}^{(n)}\rVert_{\mathcal{H}}^{p}) & \leq\mathfrak{K}_{p}+\mathfrak{K}_{p}\int_{0}^{t}\int_{[0,\tau]}\mathbb{E}(\rVert DX_{s-r}^{(n-1)}\rVert_{\mathcal{H}}^{p})\rvert\eta\rvert(\mathrm{d}r)\mathrm{d}s\\
			& \leq\sum_{k=0}^{n}\mathfrak{K}_{p}^{k+1}(\rvert\eta\rvert([0,\tau])))^{k}\frac{t^{k}}{k!}\leq\mathfrak{K}_{p}e^{\mathfrak{K}_{p}\rvert\eta\rvert([0,\tau]))T}.
		\end{align*}
		The conclusion of this proposition can now be obtained as application
		of Lemma 1.5.3 in \cite{Nualart06} and a similar reasoning used
		ins (\ref{eq:indetification_Deriv}).\end{proof}
	
	Next, we find a representation of $DX_{t}$ in terms of the kernel
	$K$ defined in (\ref{eq:kernel_fBm}). We need the following lemma.
	
	\begin{lemma} \label{exist_unique_LinearDDE} Let $(b_{t})_{t\geq0}$
		and $(N_{t})_{t\geq0}$ be two measurable processes. Assume that $b$
		is bounded and that $N$ is left (right) continuous with $\mathbb{P}(N_{t}<\infty)=1,\,\,t\geq0.$
		Then the semi-linear SDDE
		\begin{equation}
			Y_{t}=\begin{cases}
				\int_{0}^{t}\int_{[0,s]}b_{s-r}Y_{s-r}\eta(\mathrm{d}r)\mathrm{d}s+N_{t}, & t>0;\\
				0 & -\tau\leq t\leq0,
			\end{cases}\label{eq:semi-linear-SDDE}
		\end{equation}
		has at most one (up to indistinguishibility) solution. If in addition
		\begin{equation}
			\mathbb{P}(\int_{0}^{t}\rvert N_{s}\rvert\mathrm{d}s<\infty)=1,\,\,\forall\,t\geq0,\label{eq:int_semi_linear_SDDE}
		\end{equation}
		then the unique solution is given by 
		\[
		U(N)_{t}=N_{t}-\int_{0}^{t}N_{s}\frac{\partial}{\partial s}R(t,s)\mathrm{d}s,\,\,\,t\geq0,
		\]
		where $R$ is the differential resolvent of 
		\[
		\kappa(t,A,\omega)=-\int_{[0,t]}b_{t-r}(\omega)\mathbf{1}_{A}(t-r)\eta(\mathrm{d}r),\,\,t\geq0,\,\,A\in\mathcal{B}([0,+\infty)),
		\]
	\end{lemma}
	
	\begin{proof} Since $N_{t}<\infty$ a.s., the difference between
		two arbitrary solutions of (\ref{eq:semi-linear-SDDE}), say $\mathcal{E}$,
		has a version that is absolutely continuous and satisfies
		\[
		\mathcal{E}_{t}=\begin{cases}
			\int_{0}^{t}\int_{[0,s]}b_{s-r}\mathcal{E}_{s-r}\eta(\mathrm{d}r)\mathrm{d}s, & t>0;\\
			0 & -\tau\leq t\leq0.
		\end{cases}
		\]
		Hence, by Theorem 10.3.9 \cite{GripenbergLondenStaffans90} and the
		left (right) continuity of $N$, $\mathcal{E}$ is indistinguishable
		from the $0$ process. Now suppose that (\ref{eq:int_semi_linear_SDDE})
		holds and set 
		\[
		G(t,u,y):=\int_{u+y}^{t}\int_{[0,s-(y+u)]}b_{s-r}R(s-r,u+y)\eta(\mathrm{d}r)\mathrm{d}s.
		\]
		From Proposition \ref{propmeasurableresolv}, $G$ is measurable and
		bounded. Furthermore, by Fubini's Theorem and the properties of the
		differential resolvent, we get that
		\begin{align*}
			\int_{0}^{t}\int_{[0,s]}U(N)_{s-r}b_{s-r}\eta(\mathrm{d}r)\mathrm{d}s & =\int_{0}^{t}\int_{[0,t-u]}N_{u}b_{u}G(t,u,y)\eta(\mathrm{d}y)\mathrm{d}u\\
			& +\int_{0}^{t}\int_{[0,s]}N_{s-r}b_{s-r}\eta(\mathrm{d}r)\mathrm{d}s\\
			= & \int_{0}^{t}\int_{[0,t-u]}N_{u}b_{u}\left[R(t,u+y)-1\right]\eta(\mathrm{d}y)\mathrm{d}u\\
			& +\int_{0}^{t}\int_{[0,t-s]}N_{s}b_{s}\eta(\mathrm{d}r)\mathrm{d}s\\
			= & U(N)_{t}-N_{t},
		\end{align*}
		where in the lasr relation we further applied (\ref{eq:derivative_diif_resolv}).
		This concludes the proof.\end{proof}
	
	\begin{corollary}\label{Malliavin_SDDE_exact}Let the assumptions
		of Theorem \ref{prop_Malliavinder} hold. Then, up to indistinguishibility,
		\[
		D_{u}X_{t}=K(t,u)-\int_{(u)^{+}}^{t}K(s,u)\frac{\partial}{\partial s}R(t,s)\mathrm{d}s,\,\,\,u<t,
		\]
		where $K$ as in (\ref{eq:kernel_fBm}) and $R$ is the differential
		resolvent of 
		\[
		\kappa(t,A)=-\int_{[0,t]}b^{\prime}(X_{t-r})\mathbf{1}_{A}(t-r)\eta(\mathrm{d}r).
		\]
		In particular, if $T>t\geq u\geq0$, then
		\begin{equation}
			\rvert D_{u}X_{t}-(t-u)_{+}^{H-1/2}\rvert\leq C(t-u)^{H+1/2},\label{eq:MalliainDer_Estimate.}
		\end{equation}
		while for $T\geq t>0>u$, it holds
		\begin{equation}
			\rvert D_{u}X_{t}\rvert\lesssim\rvert u\rvert^{H-1/2}\mathbf{1}_{-1<u<0}+\rvert u\rvert^{H-3/2}\mathbf{1}_{u\leq-1}.\label{eq:MalliainDer_Estimate_negu}
		\end{equation}
	\end{corollary}
	
	\subsection{Exact representation of the error and a basic decomposition\label{subsec:Exact-representation-of} }
	
	In this part we find an explicit representation of the error process
	using differential resolvents (see Subsection \ref{subsec:Differential-resolvents}).
	Assume that $b$ is of class $C^{2}$ with $b^{\prime}\in C_{b}^{1}$.
	By the Mean-Value Theorem, the process $U^{n}$ satisfies the semi-linear
	delay equation 
	\[
	U_{t}^{n}=\begin{cases}
		\int_{0}^{t}\int_{[0,s]}\psi_{s-r}^{n}U_{s-r}^{n}\eta(\mathrm{d}r)\mathrm{d}s+N_{t}^{n}, & t\geq0;\\
		0, & -\tau\leq t<0,
	\end{cases}
	\]
	where 
	\[
	\psi_{s-r}^{n}=\int_{0}^{1}b^{\prime}((1-y)X_{s-r}^{n}+yX_{s-r})\mathrm{d}y,
	\]
	and 
	\begin{equation}
		N_{t}^{n}:=\int_{0}^{t}\int_{[0,\tau]}[b(X_{s-r}^{n})-b(X_{\mathcal{T}(s)-\mathcal{T}(r)}^{n})]\eta(\mathrm{d}r)\mathrm{d}s,\,\,\,t\geq0.\label{eq:Noise_def}
	\end{equation}
	Since $\psi^{n}$ is measurable and bounded, the differential resolvent
	of 
	\begin{equation}
		\kappa_{n}(t,A)=-\int_{[0,t]}\psi_{s-r}^{n}\mathbf{1}_{A}(t-r)\eta(\mathrm{d}r),\label{eq:kappa_n}
	\end{equation}
	exists and will be denoted by $R_{n}$. By Proposition \ref{propmeasurableresolv}
	we know that $R_{n}$ is jointly measurable and from \cite{GripenbergLondenStaffans90},
	Theorem 10.3.9, the error process can be written uniquely (up to indistinguishibility)
	as 
	\begin{equation}
		U_{t}^{n}=\int_{0}^{t}R_{n}(t,s)\mathrm{d}N_{s}^{n}.\label{eq:exactreperror}
	\end{equation}
	Furthermore, from Lemma \ref{exist_unique_LinearDDE} and (\ref{eq:measureinducedbyResolvent}),
	we deduce that $U$ also admits the representation
	\begin{equation}
		\begin{aligned}U_{t}^{n}= & N_{t}^{n}-\int_{0}^{t}N_{s}^{n}\frac{\partial}{\partial s}R_{n}(t,s)\mathrm{d}s\\
			= & N_{t}^{n}+\int_{0}^{t}N_{s}^{n}\int_{[0,t-s]}\psi_{s}^{n}R_{n}(t,s+r)\eta(\mathrm{d}r)\mathrm{d}s.
		\end{aligned}
		\label{eq:errordecomp2}
	\end{equation}
	Now, set 
	\begin{equation}
		S_{n}(t):=\int_{0}^{t}[b(X_{s}^{n})-b(X_{\mathcal{T}(s)}^{n})]\mathrm{d}s,\,\,t\geq0,\label{eq:def_Sn_main}
	\end{equation}
	and decompose
	\begin{equation}
		N_{t}^{n}=A_{t}^{n}+Y_{t}^{n}+D_{t}^{n},\label{eq:Decomp_noise}
	\end{equation}
	where 
	\begin{equation}
		Y_{t}^{n}=\int_{[0,t]}S_{n}(t-r)\eta(\mathrm{d}r),\label{eq:Yn_dec}
	\end{equation}
	and 
	\begin{align*}
		A_{t}^{n}:= & \int_{[0,\tau]}\int_{-\mathcal{T}(r)}^{(t-\mathcal{T}(r))\land0}[b(x_{0}(s))-b(x_{0}(\mathcal{T}(s)))]\mathrm{d}s\eta(\mathrm{d}r)\\
		& +\int_{[0,t]}\int_{t-r}^{t-\mathcal{T}(r)}[b(X_{s}^{n})-b(X_{\mathcal{T}(s)}^{n})]\mathrm{d}s\eta(\mathrm{d}r)\\
		& +\int_{(t,\left\lceil t/\Delta_{n}\right\rceil \Delta_{n})}\int_{0}^{t-\mathcal{T}(t)}[b(X_{s}^{n})-b(X_{\mathcal{T}(s)}^{n})]\mathrm{d}s\eta(\mathrm{d}r),
	\end{align*}
	as well as
	\begin{align*}
		D_{t}^{n}= & \int_{0}^{t}\int_{[0,\tau]}[b(X_{s-r}^{n})-b(X_{s-\mathcal{T}(r)}^{n})]\eta(\mathrm{d}r)\mathrm{d}s\\
		= & \int_{[0,\tau]}\left[\int_{-r}^{-\mathcal{T}(r)}b(X_{s}^{n})\mathrm{d}s-\int_{t-r}^{t-\mathcal{T}(r)}b(X_{s}^{n})\mathrm{d}s\right]\eta(\mathrm{d}r).
	\end{align*}
	It should be clear to the reader that under our assumptions (see Section
	\ref{sec:The-Euler-Maruyama-method})
	\begin{equation}
		A_{t}^{n}=\mathscr{O}_{p}^{u}(\Delta_{n});\,\,D_{t}^{n}=\mathscr{O}_{p}^{u}(\Delta_{n}).\label{eq:rate_A_Bias}
	\end{equation}
	Hence, the leading term is $Y^{n}$ whose assymptotic behaviour is
	fully described by $S_{n}$. We proceed to study this process in the
	the next section.
	
	\subsection{The core term $S_{n}$}
	
	In this part we study the asymptotic behaviour of the sequence of
	processes 
	\[
	S_{n}(t):=\int_{0}^{t}[b(X_{s}^{n})-b(X_{\mathcal{T}(s)}^{n})]\mathrm{d}s,\,\,t\geq0.
	\]
	Below, we will often use the following estimate 
	\begin{equation}
		\int_{a}^{b}\vert(x-u)^{q}-(y-u)^{q}\vert^{p}\mathrm{d}u\leq C_{p,q}(x-y)^{pq+1},\,\,x>y\geq b\geq a,\label{eq:estimateintegralbetadiff}
	\end{equation}
	where $C_{p,q}=\int_{0}^{\infty}\rvert(1+z)^{q}-z^{q}\rvert^{p}\mathrm{d}z$,
	valid for all pairs $p,q$ satisfying that $p>0$ and $0<q+\frac{1}{p}<1$.
	The following functions, defined for $u,s\geq0$ will play a fundamental
	role in our analysis
	\begin{equation}
		\begin{aligned}\varphi_{i}(s,u):= & (s-u)_{+}^{\beta}-(t_{i-1}-u)_{+}^{\beta};\\
			\psi_{i}^{n}(u):= & \int_{t_{i-1}}^{t_{i}}\varphi_{i}(s,u)\mathrm{d}s;\\
			\chi_{i}^{n}(u):= & (t_{i-1}-u)_{+}^{\beta}\psi_{i}^{n}(u)+\frac{1}{2}\int_{t_{i-1}}^{t_{i}}\varphi_{i}(s,u)^{2}\mathrm{d}s;\\
			\gamma_{i}^{n}(u):= & \beta\int_{t_{i-1}}^{t_{i}}(y-u)_{+}^{\beta-1}[(t_{i}-y)-(y-t_{i-1})]\mathrm{d}y,\,\,\,H>1/2,
		\end{aligned}
		\label{eq:def_rough_function-1}
	\end{equation}
	where we have let $t_{i}=i\Delta_{n},$ $i=0,1,\ldots$ and set $\beta:=H-1/2$.
	We start by analyzing $S_{n}$ when $H\geq1/2$. 
	
	\begin{theorem} \label{thmkey}Let $S_{n}$ be as in (\ref{eq:def_Sn_main}).
		Assume that $H\geq1/2$ and that $b$ is of class $\mathcal{C}^{2}$
		with $b^{\prime}\in\mathcal{C}_{b}^{1}$. Then, for any $p\geq1$
		and $T>0$
		\[
		\sup_{0\leq t\leq T}\mathbb{E}(\rvert S_{n}(t)\rvert^{p})^{1/p}\leq C\Delta_{n}.
		\]
		Furthermore,
		\begin{enumerate}
			\item If $H>1/2$, then 
			\[
			\frac{1}{\Delta_{n}}S_{n}\overset{u.c.p}{\longrightarrow}\frac{1}{2}[b(X_{\cdot})-b(X_{0})].
			\]
			\item If $H=1/2$, then 
			\begin{equation}
				\frac{1}{\Delta_{n}}S_{n}\overset{\mathcal{F}-C[0,T]}{\Longrightarrow}\frac{1}{2}[b(X_{\cdot})-b(X_{0})]+\frac{1}{\sqrt{12}}\int_{0}^{\cdot}b^{\prime}(X_{s})\mathrm{d}\tilde{W}_{s},\label{eq:conv_Sn_BM}
			\end{equation}
			where $\tilde{W}$ as in Theorem \ref{thmerrordisGaussiannonrough}.
		\end{enumerate}
	\end{theorem}
	
	\begin{proof}Plainly 
		\begin{equation}
			S_{n}(t):=\int_{0}^{t}b^{\prime}(X_{\mathcal{T}(s)}^{n})[X_{s}^{n}-X_{\mathcal{T}(s)}^{n}]\mathrm{d}s+\frac{1}{2}\int_{0}^{t}b^{\prime\prime}(\theta^{n}(s)X_{s}^{n}+(1-\theta^{n}(s))X_{\mathcal{T}(s)}^{n})(X_{s}^{n}-X_{\mathcal{T}(s)}^{n})^{2}\mathrm{d}s,\label{eq:first_order_Sn}
		\end{equation}
		for some random number $\theta^{n}(s)\in[0,1]$.
		
		\textbf{Assume that $H>1/2$:} In this situation, (\ref{eq:naiveineq})
		and Lemma \ref{estimatesincXn} in combination with \eqref{eq:first_order_Sn} imply that 
		\begin{equation}
			\begin{aligned}S_{n}(t) & :=\int_{0}^{t}b^{\prime}(X_{\mathcal{T}(s)})[B_{s}-B_{\mathcal{T}(s)}]\mathrm{d}s+\mathscr{O}_{p}^{u}(\Delta_{n}^{2H})\\
				& +\int_{0}^{t}b^{\prime}(X_{\mathcal{T}(s)})\int_{\mathcal{T}(s)}^{s}\int_{[0,\tau]}b(X_{\mathcal{T}(u)-\mathcal{T}(y)}^{n})\eta(\mathrm{d}y)\mathrm{d}u\mathrm{d}s.
			\end{aligned}
			\label{eq:firstoderdecom_Sn}
		\end{equation}
		Remark \ref{rmk_momentEuler} indicates that the last summand is $\mathscr{O}_{p}^{u}(\Delta_{n}).$ Furthermore,
		by arguing as below (see equation (\ref{eq:limitA}) and the subsequent
		reasoning), we easily see that, after normalizing by $\frac{1}{\Delta_{n}}$,
		such term converges pointwise in probability towards the process
		\[
		\frac{1}{2}\int_{0}^{t}b^{\prime}(X_{s})\int_{[0,\tau]}b(X_{s-y})\eta(\mathrm{d}y)\mathrm{d}s,\,\,\,\forall\,t\geq0.
		\]
		Note that this convergence is also uniformly in compacts
		due to Dini's Theorem. Now, by Corollary 7.2 in \cite{HuLiuNualart16}
		we further have that
		\[
		\frac{1}{\Delta_{n}}\int_{0}^{t}b^{\prime}(X_{\mathcal{T}(s)})[B_{s}-B_{\mathcal{T}(s)}]\mathrm{d}s\overset{L^{p}}{\rightarrow}\frac{1}{2}\int_{0}^{t}b^{\prime}(X_{s})\mathrm{d}B_{s},\,\,\,\forall\,t\geq0.
		\]
		Since the limit is continuous, according to Lemma A2.1 in \cite{DonKutz96}
		and Corollary 2.2 in \cite{BurdzySwanson10}, we only need to show
		that for every $T>0$, and $\mathsf{I},\mathsf{J}=0,1,\ldots,[T/\Delta_{n}]$,
		with $\mathsf{I}<\mathsf{J}$, it holds that 
		\begin{equation}
			\mathbb{E}(\rvert\int_{\mathsf{I}\Delta_{n}}^{\mathsf{J}\Delta_{n}}b^{\prime}(X_{\mathcal{T}(s)})[B_{s}-B_{\mathcal{T}(s)}]\mathrm{d}s\rvert^{2})^{1/2}\leq\Delta_{n}C[(\mathsf{J}-\mathsf{I})\Delta_{n}]^{H}.\label{eq:tightness_nonrough}
		\end{equation}
		Using (\ref{eq:decompFLSM}) and the Stochastic Fubini Theorem (see
		for instance Theorem 3.1 in \cite{BasseBN11}), we write 
		\begin{equation}
			\begin{aligned}\int_{\mathsf{I}\Delta_{n}}^{\mathsf{J}\Delta_{n}}b^{\prime}(X_{\mathcal{T}(s)})[B_{s}-B_{\mathcal{T}(s)}]\mathrm{d}s= & \sum_{i=\mathsf{I}+1}^{\mathsf{J}}b^{\prime}(X_{t_{i-1}})\int_{t_{i-1}}^{t_{i}}v_{u}(t_{i}-u)\mathrm{d}u\\
				& +\sum_{i=\mathsf{I}+1}^{\mathsf{J}}b^{\prime}(X_{t_{i-1}})\int_{0}^{t_{i}}\psi_{i}^{n}(u)\mathrm{d}W_{u}\\
				=: & \mathcal{I}_{\mathsf{I},\mathsf{J}}^{1}+\mathcal{I}_{\mathsf{I},\mathsf{J}}^{2}.
			\end{aligned}
			\label{eq:basic_decomposition_Sn}
		\end{equation}
		Obviously 
		\[
		\rvert\mathcal{I}_{\mathsf{I},\mathsf{J}}^{1}\rvert\leq\Delta_{n}\int_{\mathsf{I}\Delta_{n}}^{\mathsf{J}\Delta_{n}}\rvert v_{u}\rvert\mathrm{d}u.
		\]
		In view that $v$ is Gaussian, we further have that $\mathbb{E}(\vert v_{t}\vert^{p})=C_{H,p}t^{p(H-1)}.$
		Hence, for all $H\in(0,1)$ and every $p\geq1$
		\begin{equation}
			\mathbb{E}(\rvert\mathcal{I}_{\mathsf{I},\mathsf{J}}^{1}\rvert^{p})^{1/p}\leq\Delta_{n}C\int_{\mathsf{I}\Delta_{n}}^{\mathsf{J}\Delta_{n}}u^{H-1}\mathrm{d}u\leq\Delta_{n}C[(\mathsf{J}-\mathsf{I})\Delta_{n}]^{H}.\label{eq:BV_tightness}
		\end{equation}
		Next, we write
		\begin{equation}
			\begin{aligned}\mathcal{I}_{\mathsf{I},\mathsf{J}}^{2}= & \sum_{i=\mathsf{I}+1}^{\mathsf{J}}\int_{t_{i-1}}^{t_{i}}b^{\prime}(X_{t_{i-1}})\psi_{i}^{n}(u)\mathrm{d}W_{u}+\sum_{i=\mathsf{I}+1}^{\mathsf{J}}\int_{0}^{t_{i-1}}b^{\prime}(X_{t_{i-1}})\psi_{i}^{n}(u)\delta W_{u}\\
				& +\sum_{i=\mathsf{I}+1}^{\mathsf{J}}\int_{0}^{t_{i-1}}b^{\prime\prime}(X_{t_{i-1}})D_{u}X_{t_{i-1}}\psi_{i}^{n}(u)\mathrm{d}u\\
				=: & \mathcal{I}_{\mathsf{I},\mathsf{J}}^{2,1}+\mathcal{I}_{\mathsf{I},\mathsf{J}}^{2,2}+\mathcal{I}_{\mathsf{I},\mathsf{J}}^{2,3}.
			\end{aligned}
			\label{eq:decomp_I2}
		\end{equation}
		Using that for $t_{i-1}\leq u\leq t_{i}$, $\psi_{i}^{n}(u)=\frac{1}{H+1/2}(t_{i}-u)^{H+1/2}$
		and the Burkholder-Davis-Gundy inequality, we deduce that for all
		$p\geq1$
		\begin{equation}
			\mathbb{E}(\rvert\mathcal{I}_{\mathsf{I},\mathsf{J}}^{2,1}\rvert^{p})^{1/p}\leq C\Delta_{n}^{H+1/2}\left[(\mathsf{J}-\mathsf{I})\Delta_{n}\right]^{1/2}\leq\Delta_{n}C\left[(\mathsf{J}-\mathsf{I})\Delta_{n}\right]^{H},\label{eq:tighness_nonrough_I21}
		\end{equation}
		where in the last step we used the fact that $\mathsf{J}-\mathsf{I}\geq1$
		and $H>1/2$. Now we note that
		\begin{equation}
			\begin{aligned}\mathcal{I}_{\mathsf{I},\mathsf{J}}^{2,3}= & \sum_{i=\mathsf{I}+1}^{\mathsf{J}}\sum_{k=0}^{i-1}\int_{t_{k-1}}^{t_{k}}b^{\prime\prime}(X_{t_{i-1}})D_{u}X_{t_{i-1}}\psi_{i}^{n}(u)\mathrm{d}u\\
				= & \sum_{k=0}^{\mathsf{J}-1}\int_{t_{k-1}}^{t_{k}}\sum_{i=\mathsf{I}\lor k+1}^{\mathsf{J}}b^{\prime\prime}(X_{t_{i-1}})D_{u}X_{t_{i-1}}\psi_{i}^{n}(u)\mathrm{d}u
			\end{aligned}
			\label{eq:decomp_I23}
		\end{equation}
		and similarly 
		\begin{equation}
			\begin{aligned}\mathcal{I}_{\mathsf{I},\mathsf{J}}^{2,2}= & \sum_{k=0}^{\mathsf{J}-1}\int_{t_{k-1}}^{t_{k}}\sum_{i=\mathsf{I}\lor k+1}^{\mathsf{J}}b^{\prime}(X_{t_{i-1}})\psi_{i}^{n}(u)\delta W_{u}\end{aligned}
			\label{eq:decomp_I22}
		\end{equation}
		Relation (\ref{eq:MalliainDer_Estimate.}) along with part 3. of Lemma
		\ref{Lemma_estimates_phi} below imply that 
		\[
		\rvert\mathcal{I}_{\mathsf{I},\mathsf{J}}^{2,3}\rvert\lesssim\int_{0}^{(\mathsf{J}-1)\Delta_{n}}\rvert f_{n}(\mathsf{J}\Delta_{n},u)-f_{n}(\mathsf{I}\Delta_{n},u)\rvert\mathrm{d}u\lesssim\Delta_{n}\left[(\mathsf{J}-\mathsf{I})\Delta_{n}\right].
		\]
		Note now that for all $i=1,\ldots,$ the process $b^{\prime}(X_{t_{i-1}})\psi_{i}^{n}(\cdot)$
		belongs $\mathbb{D}^{1,p}(\mathcal{H})$ for all $p>1$. Hence, by
		(\ref{eq:meyers_ine}) and Corollary \ref{Malliavin_SDDE_exact} we
		infer that
		\begin{align*}
			\mathbb{E}(\rvert\mathcal{I}_{\mathsf{I},\mathsf{J}}^{2,2}\rvert^{2})\lesssim & \int_{0}^{(\mathsf{J}-1)\Delta_{n}}\rvert f_{n}(\mathsf{J}\Delta_{n},u)-f_{n}(\mathsf{I}\Delta_{n},u)\rvert^{2}\mathrm{d}u\\
			\leq & C\Delta_{n}^{2}\left[(\mathsf{J}-\mathsf{I})\Delta_{n}\right]^{2H}
		\end{align*}
		where in the last inequality we applied part 3. of Lemma \ref{Lemma_estimates_phi}.
		The previous two estimates in combination with (\ref{eq:BV_tightness})
		establish the validity of (\ref{eq:tightness_nonrough}) concluding
		this the proof of (1).
		
		\textbf{Assume that $H=1/2$:} We note first that in this case $B=W$,
		so that $X$ and $X^{n}$ are continuous $(\mathcal{F}_{t}^{W})_{t\geq0}$-semimartingales
		and
		\begin{equation}
			S_{n}(t):=\mathfrak{C}_{t}^{n,1}+\mathfrak{C}_{t}^{n,2}+M_{t}^{n},\label{eq:first_order_Sn-1}
		\end{equation}
		where 
		\begin{align*}
			\mathfrak{C}_{t}^{n,1} & :=\frac{1}{2}\int_{0}^{t}b^{\prime\prime}(\theta^{n}(s)X_{s}^{n}+(1-\theta^{n}(s))X_{\mathcal{T}(s)}^{n})(X_{s}^{n}-X_{\mathcal{T}(s)}^{n})^{2}\mathrm{d}s\\
			\mathfrak{C}_{t}^{n,2} & :=\int_{0}^{t}b^{\prime}(X_{\mathcal{T}(s)}^{n})\int_{\mathcal{T}(s)}^{s}\int_{[0,\tau]}b(X_{\mathcal{T}(u)-\mathcal{T}(y)}^{n})\eta(\mathrm{d}y)\mathrm{d}u\mathrm{d}s,
		\end{align*}
		and 
		\[
		M_{t}^{n}:=\int_{0}^{t}b^{\prime}(X_{\mathcal{T}(s)}^{n})[W_{s}-W_{\mathcal{T}(s)}]\mathrm{d}s.
		\]
		Our goal is to show that 
		\begin{align}
			1.\, & \frac{1}{\Delta_{n}}\mathfrak{C}^{n,1}\overset{u.c.p}{\longrightarrow}\frac{1}{4}\int_{0}^{\cdot}b^{\prime\prime}(X_{s})\mathrm{d}s;\label{eq:BM_drift1}\\
			2.\, & \frac{1}{\Delta_{n}}\mathfrak{C}^{n,2}\overset{u.c.p}{\longrightarrow}\frac{1}{2}\int_{0}^{t}b^{\prime}(X_{s})\int_{[0,\tau]}b(X_{s-y})\eta(\mathrm{d}y)\mathrm{d}s;\label{eq:BM_drift2}\\
			3.\, & \frac{1}{\Delta_{n}}M^{n}\overset{\mathcal{F}-\mathcal{C}([0,T])}{\Longrightarrow}\frac{1}{2}\int_{0}^{\cdot}b^{\prime}(X_{s})\mathrm{d}W_{s}+\frac{1}{\sqrt{12}}\int_{0}^{\cdot}b^{\prime}(X_{s})\mathrm{d}\tilde{W}_{s}.\label{eq:BM_mtgpart}
		\end{align}
		If this were true, we would have from (\ref{SDDE}) 
		\[
		\frac{1}{\Delta_{n}}S_{n}\overset{\mathcal{F}-\mathcal{C}([0,T])}{\Longrightarrow}\frac{1}{2}\int_{0}^{\cdot}b^{\prime}(X_{s})\mathrm{d}X_{s}+\frac{1}{4}\int_{0}^{\cdot}b^{\prime\prime}(X_{s})\mathrm{d}s+\frac{1}{\sqrt{12}}\int_{0}^{\cdot}b^{\prime}(X_{s})\mathrm{d}\tilde{W}_{s},
		\]
		which is exactly (\ref{eq:conv_Sn_BM}) due to Itô's formula. We have
		already seen that (\ref{eq:BM_drift2}) holds. By Dini's Theorem we
		only need to check that (\ref{eq:BM_drift1}) holds for each $t\geq0$.
		To that end, first observe that from (\ref{eq:naiveineq}), Lemma \ref{estimatesincXn},
		and the fact that $b^{\prime}\in\mathcal{C}_{b}^{1}$, it holds that
		\[
		\mathfrak{C}_{t}^{n,1}=\frac{1}{2}\int_{0}^{t}b^{\prime\prime}(X_{\mathcal{T}(s)})(X_{s}^{n}-X_{\mathcal{T}(s)}^{n})^{2}\mathrm{d}s+\mathrm{o}_{\mathbb{P}}(\Delta_{n}).
		\]
		A standard aplication of Itô's formula along with Lemma \ref{estimatesincXn}
		give us that 
		\begin{align*}
			\frac{1}{2}\int_{0}^{t}b^{\prime\prime}(X_{\mathcal{T}(s)})(X_{s}^{n}-X_{\mathcal{T}(s)}^{n})^{2}\mathrm{d}s= & \frac{1}{2}\int_{0}^{t}b^{\prime}(X_{\mathcal{T}(s)})(s-\mathcal{T}(s))\mathrm{d}s+\mathscr{O}_{p}^{u}(\Delta_{n}^{3/2})\\
			& +\int_{0}^{t}b^{\prime\prime}(X_{\mathcal{T}(s)})\int_{\mathcal{T}(s)}^{s}(X_{u}^{n}-X_{\mathcal{T}(s)}^{n})\mathrm{d}W_{u}\mathrm{d}s.
		\end{align*}
		Clearly
		\[
		\frac{1}{2\Delta_{n}}\int_{0}^{t}b^{\prime\prime}(X_{\mathcal{T}(s)})(s-\mathcal{T}(s))\mathrm{d}s\overset{u.c.p}{\longrightarrow}\frac{1}{4}\int_{0}^{t}b^{\prime\prime}(X_{s})\mathrm{d}s.
		\]
		On the other hand, by a further application of the stochastic Fubini
		theorem and Lemma \ref{estimatesincXn} we get
		\begin{align*}
			\mathfrak{N}_{t}^{n}:= & \int_{0}^{t}b^{\prime\prime}(X_{\mathcal{T}(s)})\int_{\mathcal{T}(s)}^{s}(X_{u}^{n}-X_{\mathcal{T}(s)}^{n})\mathrm{d}W_{u}\mathrm{d}s\\
			= & \sum_{i=1}^{[t/\Delta_{n}]}\int_{t_{i-1}}^{t_{i}}b^{\prime\prime}(X_{t_{i-1}})(X_{u}^{n}-X_{t_{i-1}}^{n})(t_{i}-u)\mathrm{d}W_{u}+\mathscr{O}_{p}^{u}(\Delta_{n}^{3/2}).
		\end{align*}
		Since $\mathbb{E}\left\{ [b^{\prime}(X_{t_{i-1}})(X_{u}^{n}-X_{t_{i-1}}^{n})(t_{i}-u)]^{2}\right\} \leq C\Delta_{n}^{3},$
		Lemma 2.2.11 in \cite{JacProt11} implies that $\mathfrak{N}_{t}^{n}=\mathrm{o}_{\mathbb{P}}(\Delta_{n})$, which completes the proof of \ref{eq:BM_drift1}).
		It is left to show that (\ref{eq:BM_mtgpart}) holds. We begin by
		checking that $M^{n}$ is tight on $\mathcal{C}([0,T])$. Specifically
		(see Corollary 16.9 in \cite{Kallenberg02}), we will establish that for
		every $p\geq1$
		\begin{equation}
			\mathbb{E}(\rvert M_{t}^{n}-M_{v}^{n}\rvert^{p})\lesssim\Delta_{n}^{p}(t-v)^{p/2},\,\,0\leq v\leq t\leq T.\label{eq:tightness_M_BM}
		\end{equation}
		Note that if $0\leq t-v<\Delta_{n}$, Jensen's inequality implies
		that 
		\[
		\mathbb{E}(\rvert M_{t}^{n}-M_{v}^{n}\rvert^{p})\leq\rVert b^{\prime}\rVert_{\infty}^{p}\Delta_{n}^{p/2}(t-v)^{p}\leq\rVert b^{\prime}\rVert_{\infty}^{p}\Delta_{n}^{p}(t-v)^{p/2}.
		\]
		If instead $t-v\geq\Delta_{n}$, by arguing as in the case of $\mathfrak{N}^{n}$,
		we have that 
		\begin{equation}
			\begin{aligned}M_{t}^{n}= & \sum_{i=1}^{[t/\Delta_{n}]}b^{\prime}(X_{t_{i-1}})\int_{t_{i-1}}^{t_{i}}(t_{i}-u)\mathrm{d}W_{u}+\int_{\mathcal{T}(t)}^{t}b^{\prime}(X_{\mathcal{T}(t)}^{n})[W_{s}-W_{\mathcal{T}(t)}]\mathrm{d}s\\
				= & \int_{0}^{t}b^{\prime}(X_{\mathcal{T}(s)}^{n})(\mathcal{V}(u)-u)\mathrm{d}W_{u}+b^{\prime}(X_{\mathcal{T}(t)}^{n})\int_{\mathcal{T}(t)}^{t}(t-u)\mathrm{d}W_{u}\\
				& -b^{\prime}(X_{\mathcal{T}(t)}^{n})\int_{\mathcal{T}(t)}^{t}(\mathcal{V}(t)-u)\mathrm{d}W_{u},
			\end{aligned}
			\label{eq:dec_mtg_BM}
		\end{equation}
		where $\mathcal{V}(s)=\lceil s/\Delta_{n}\rceil\Delta_{n}$. Hence,
		by the Burkholder-Davis-Gundy inequality, we deduce that 
		\[
		\mathbb{E}(\rvert M_{t}^{n}-M_{v}^{n}\rvert^{p})\leq C_{p}\rVert b^{\prime}\rVert_{\infty}^{p}\left\{ \Delta_{n}^{p}(t-v)^{p/2}+\Delta_{n}^{\frac{3}{2}p}\right\} \leq C_{p}\rVert b^{\prime}\rVert_{\infty}^{p}\Delta_{n}^{p}(t-v)^{p/2},
		\]
		where in the last inequality we also used that $t-v\geq\Delta_{n}$.
		Finally, in view of the decomposition (\ref{eq:dec_mtg_BM}) and since
		\begin{align*}
			1.\, & \frac{1}{\Delta_{n}^{2}}\sum_{i=1}^{[t/\Delta_{n}]}b^{\prime}(X_{t_{i-1}})^{2}\int_{t_{i-1}}^{t_{i}}(t_{i}-u)^{2}\mathrm{d}u\overset{\mathbb{P}}{\longrightarrow}\frac{1}{3}\int_{0}^{t}b^{\prime}(X_{s})^{2}\mathrm{d}s;\\
			2.\, & \frac{1}{\Delta_{n}}\sum_{i=1}^{[t/\Delta_{n}]}b^{\prime}(X_{t_{i-1}})\int_{t_{i-1}}^{t_{i}}(t_{i}-u)\mathrm{d}u\overset{\mathbb{P}}{\longrightarrow}\frac{1}{2}\int_{0}^{t}b^{\prime}(X_{s})\mathrm{d}s;\\
			3.\, & \frac{1}{\Delta_{n}^{p}}\sum_{i=1}^{[t/\Delta_{n}]}b^{\prime}(X_{t_{i-1}})^{p}\left(\int_{t_{i-1}}^{t_{i}}(t_{i}-u)^{2}\mathrm{d}u\right)^{p/2}\overset{\mathbb{P}}{\longrightarrow}0,\,\,\forall\,p>2;
		\end{align*}
		we conclude from Theorem IX.7.28 in \cite{JacShri02} (with $Z=W$)
		and Proposition 3.20 in \cite{HauslerLuschgy15} that (\ref{eq:BM_mtgpart})
		is indeed true.\end{proof} 
	
	Next, we consider the rough case, i.e. $0<H<1/2$. We introduce the
	following processes
	\begin{equation}
		\begin{aligned}L_{t}^{n}:= & \int_{0}^{t}b^{\prime}(X_{\mathcal{T}(s)})(Z_{s}-Z_{\mathcal{T}(s)})\mathrm{d}s;\\
			Q_{t}^{n}:= & \int_{0}^{t}b^{\prime\prime}(X_{\mathcal{T}(s)})(Z_{s}-Z_{\mathcal{T}(s)})^{2}\mathrm{d}s,
		\end{aligned}
		\label{eq:def_L_Q}
	\end{equation}
	and set 
	\begin{equation}
		S_{n}^{\star}(t):=L_{t}^{n}+\frac{1}{2}Q_{t}^{n}.\label{eq:def_Sn_star}
	\end{equation}
	Our next goal is to show that for any $0<H<1/2$, $S_{n}^{\star}(t)=\mathscr{O}_{p}^{u}(\Delta_{n}^{(H+1/2)\land3H})$.
	More precisely,
	
	\begin{theorem} \label{thmkey-1}Let $S_{n}$ be as in (\ref{eq:def_Sn_star})
		and assume that $0<H<1/2$. Suppose in addition that $b$ is of class
		$\mathcal{C}^{3}$ with $b^{\prime}\in\mathcal{C}_{b}^{2}$. Then,
		for every $T>0$ and $p\geq1$
		\[
		\mathbb{E}(\rvert S_{n}^{\star}(t)-S_{n}^{\star}(v)\rvert^{p})^{1/p}\lesssim\Delta_{n}^{(H+1/2)\land3H}(t-v)^{1/2},\,\,\,0\leq v\leq t\leq T.
		\]
		
	\end{theorem}
	\begin{proof}Let us start showing that for any $T>0$, $p\geq1$, and $\mathsf{I},\mathsf{J}=0,1,\ldots,[T/\Delta_{n}]$,
		with $\mathsf{I}<\mathsf{J}$, it holds that 
		\begin{equation}
			\mathbb{E}(\rvert S_{n}^{\star}(\mathsf{J}\Delta_{n})-S_{n}^{\star}(\mathsf{I}\Delta_{n})\rvert^{p})\leq C\Delta_{n}^{(H+1/2)\land3H}[\Delta_{n}(\mathrm{J}-\mathrm{I})]{}^{p/2}.\label{eq:in_rationals}
		\end{equation}
		By Itô's formula, for all $t_{i-1}\leq s<t_{i}$, $i=1,2,\ldots$,
		we have that 
		\[
		(Z_{s}-Z_{t_{i-1}})^{2}=2\int_{0}^{t_{i}}\int_{0}^{u}\varphi_{i}(s,u)\varphi_{i}(s,v)\mathrm{d}W_{v}\mathrm{d}W_{u}+\int_{0}^{t_{i}}\varphi_{i}(s,u)^{2}\mathrm{d}u.
		\]
		Hence, by letting 
		\[
		\Psi_{i}^{n}(u):=\int_{t_{i-1}}^{t_{i}}\varphi_{i}(s,u)\int_{0}^{u}\varphi_{i}(s,v)\mathrm{d}W_{v}\mathrm{d}s,
		\]
		and arguing as in (\ref{eq:decomp_I2}), we obtain (with the same
		notation) that 
		\[
		S_{n}^{\star}(\mathsf{J}\Delta_{n})-S_{n}^{\star}(\mathsf{I}\Delta_{n})=\mathcal{I}_{\mathsf{I},\mathsf{J}}^{2,1}+\mathcal{I}_{\mathsf{I},\mathsf{J}}^{2,2}+\mathcal{I}_{\mathsf{I},\mathsf{J}}^{2,3}+\mathcal{I}_{\mathsf{I},\mathsf{J}}^{3,1}+\mathcal{I}_{\mathsf{I},\mathsf{J}}^{3,2}+\mathcal{I}_{\mathsf{I},\mathsf{J}}^{3,3}+\mathcal{I}_{\mathsf{I},\mathsf{J}}^{3,4};
		\]
		in which 
		\begin{align*}
			\mathcal{I}_{\mathsf{I},\mathsf{J}}^{3,1} & =\sum_{i=\mathsf{I}+1}^{\mathsf{J}}\int_{t_{i-1}}^{t_{i}}b^{\prime\prime}(X_{t_{i-1}})\Psi_{i}^{n}(u)\mathrm{d}W_{u};\\
			\mathcal{I}_{\mathsf{I},\mathsf{J}}^{3,2} & =\sum_{i=\mathsf{I}+1}^{\mathsf{J}}\int_{0}^{t_{i-1}}b^{\prime\prime}(X_{t_{i-1}})\Psi_{i}^{n}(u)\delta W_{u};\\
			\mathcal{I}_{\mathsf{I},\mathsf{J}}^{3,3} & =\sum_{i=\mathsf{I}+1}^{\mathsf{J}}\int_{0}^{t_{i-1}}b^{\prime\prime\prime}(X_{t_{i-1}})D_{u}X_{t_{i-1}}\Psi_{i}^{n}(u)\mathrm{d}u;\\
			\mathcal{I}_{\mathsf{I},\mathsf{J}}^{3,4} & =\frac{1}{2}\sum_{i=\mathsf{I}+1}^{\mathsf{J}}\int_{0}^{t_{i}}b^{\prime\prime}(X_{t_{i-1}})\int_{t_{i-1}}^{t_{i}}\varphi_{i}(s,u)^{2}\mathrm{d}s\mathrm{d}u.
		\end{align*}
		Note that the first part of (\ref{eq:tighness_nonrough_I21}) remains
		valid, i.e. 
		\begin{equation}
			\mathbb{E}(\rvert\mathcal{I}_{\mathsf{I},\mathsf{J}}^{2,1}\rvert^{p})^{1/p}\leq C\Delta_{n}^{H+1/2}\left[\Delta_{n}(\mathsf{J}-\mathsf{I})\right]^{1/2}.\label{eq:I21_rough_lpest}
		\end{equation}
		Furthermore, by part 4 of Lemma \ref{Lemma_estimates_phi} below and
		Corollary \ref{Malliavin_SDDE_exact} we infer that
		\[
		\rvert\mathcal{I}_{\mathsf{I},\mathsf{J}}^{2,3}+\mathcal{I}_{\mathsf{I},\mathsf{J}}^{3,4}\rvert\lesssim\Delta_{n}^{H+1/2}\left\{ [\Delta_{n}(\mathrm{J}-\mathrm{I})]^{H+1/2}+[\Delta_{n}(\mathrm{J}-\mathrm{I})]\right\} .
		\]
		Now set
		\[
		\mathbb{I}_{\mathsf{I},\mathsf{J}}^{2,2}:=\sum_{l=1}^{\mathsf{J}-1}\sum_{k=1}^{\mathsf{J}-1}\int_{t_{l-1}}^{t_{l}}\int_{t_{k-1}}^{t_{k}}\left(\sum_{i=\mathsf{I}\lor k+1}^{\mathsf{J}}(t_{i-1}-v)^{\beta}\rvert\psi_{i}^{n}(u)\rvert\mathbf{1}_{i\geq l+1}\right)^{2}\mathrm{d}u\mathrm{d}v,
		\]
		Then, 
		\begin{equation}
			\begin{aligned}\mathbb{E}(\rvert\mathcal{I}_{\mathsf{I},\mathsf{J}}^{2,2}\rvert^{p}) & \lesssim\left(\int_{0}^{(\mathsf{J}-1)\Delta_{n}}\rvert f_{n}(\mathsf{J}\Delta_{n},u)-f_{n}(\mathsf{I}\Delta_{n},u)\rvert^{2}\mathrm{d}u\right)^{p/2}+(\mathbb{I}_{\mathsf{I},\mathsf{J}}^{2,2})^{p/2}\\
				& \lesssim\Delta_{n}^{p(H+1/2)}\left[\Delta_{n}(\mathsf{J}-\mathsf{I})\right]^{p/2}+(\mathbb{I}_{\mathsf{I},\mathsf{J}}^{2,2})^{p/2},
			\end{aligned}
			\label{eq:appl_meyersine_rough}
		\end{equation}
		due to (\ref{eq:meyers_ine}), Corollary \ref{Malliavin_SDDE_exact},
		and part 3 of Lemma \ref{Lemma_estimates_phi}. By expanding the squared
		and changing the order of summation we conclude that
		\begin{equation}
			\begin{aligned}\mathbb{I}_{\mathsf{I},\mathsf{J}}^{2,2} & \leq C\sum_{k=1}^{\mathsf{J}-1}\sum_{i,i^{\prime}=\mathsf{I}\lor k+1}^{\mathsf{J}}\int_{t_{k-1}}^{t_{k}}\rvert\psi_{i}^{n}(u)\rvert\rvert\psi_{i^{\prime}}^{n}(u)\rvert\sum_{l=1}^{i\land i^{\prime}}\int_{0}^{t_{i-1}\land t_{i^{\prime}-1}}(t_{i-1}\land t_{i^{\prime}-1}-v)^{2\beta}\mathrm{d}v\mathrm{d}u\\
				& \leq CT^{2H}\int_{0}^{(\mathsf{J}-1)\Delta_{n}}\rvert f_{n}(\mathsf{J}\Delta_{n},u)-f_{n}(\mathsf{I}\Delta_{n},u)\rvert^{2}\mathrm{d}u\\
				& \lesssim\Delta_{n}^{H+1/2}\left[\Delta_{n}(\mathsf{J}-\mathsf{I})\right],
			\end{aligned}
			\label{eq:norm_bound_meyet}
		\end{equation}
		where in the last step we applied once again part 3 of Lemma \ref{Lemma_estimates_phi}.
		To deal with the rest of the terms we will use that for any $p\geq1$
		and $u\leq t_{i-1}$
		\begin{equation}
			\mathbb{E}(\rvert\Psi_{i}^{n}(u)\rvert^{p})^{1/p}\leq C_{p}\int_{t_{i-1}}^{t_{i}}(\int_{0}^{t_{i-1}}\varphi_{i}(s,v)^{2}\mathrm{d}v)^{1/2}\rvert\varphi_{i}(s,u)\rvert\mathrm{d}s\leq  C_{p}\Delta_{n}^{H}\rvert\psi_{i}^{n}(u)\rvert,
			\label{eq:Lp_estimate_PHI}
		\end{equation}
		which can be obtained by applying Minskowski's inequality, the Burkholder-Davis-Gundy
		inequality, and (\ref{eq:estimateintegralbetadiff}). This bound along
		with Corollary \ref{Malliavin_SDDE_exact} and part 4 of Lemma \ref{Lemma_estimates_phi}
		result in 
		\begin{align*}
			\mathbb{E}(\rvert\mathcal{I}_{\mathsf{I},\mathsf{J}}^{3,3}\rvert^{p})^{1/p}\lesssim & \Delta_{n}^{H}\sum_{i=\mathsf{I}+1}^{\mathsf{J}}\int_{0}^{t_{i-1}}(t_{i-1}-u)^{\beta}\rvert\psi_{i}^{n}(u)\rvert\mathrm{d}u+\Delta_{n}^{H}\sum_{i=\mathsf{I}+1}^{\mathsf{J}}\int_{0}^{t_{i-1}}\rvert\psi_{i}^{n}(u)\rvert\mathrm{d}u\\
			\lesssim & \Delta_{n}^{3H}\left[\Delta_{n}(\mathsf{J}-\mathsf{I})\right]+\Delta_{n}^{2H+1/2}\left[\Delta_{n}(\mathsf{J}-\mathsf{I})\right],
		\end{align*}
		as well as that 
		\[
		\mathbb{E}(\rvert\mathcal{I}_{\mathsf{I},\mathsf{J}}^{3,1}\rvert^{p})\lesssim\mathbb{E}\left[\left(\sum_{i=\mathsf{I}+1}^{\mathsf{J}}\int_{t_{i-1}}^{t_{i}}\Psi_{i}^{n}(u)^{2}\mathrm{d}u\right)^{p/2}\right]\lesssim\Delta_{n}^{p(2H+1/2)}\left[\Delta_{n}(\mathsf{J}-\mathsf{I})\right]^{p/2}.
		\]
		Finally, (\ref{eq:Lp_estimate_PHI}) allows us to substitute $\psi_{i}^{n}(u)$ by $\Psi_{i}^{n}(u)$
		in (\ref{eq:appl_meyersine_rough}) and (\ref{eq:norm_bound_meyet}), so that
		
		\[
		\mathbb{E}(\rvert\mathcal{I}_{\mathsf{I},\mathsf{J}}^{3,2}\rvert^{p})\lesssim\lesssim\Delta_{n}^{p(2H+1/2)}\left[\Delta_{n}(\mathsf{J}-\mathsf{I})\right]^{p/2},
		\]
		which concludes our argument for \eqref{eq:in_rationals}. To end the proof, let us now take arbitrary $0\leq v<t\leq T.$ If
		$0<t-v<\Delta_{n}$, we easily have
		\[
		\mathbb{E}(\rvert S_{n}^{\star}(t)-S_{n}^{\star}(v)\rvert^{p})^{1/p}\lesssim\Delta_{n}^{H}(t-v)\leq\Delta_{n}^{H+1/2}(t-v)^{1/2},
		\]
		thanks to \eqref{eq:estimateintegralbetadiff}.	If $t-v\geq\Delta_{n}$, (\ref{eq:in_rationals}) now implies that 
		\begin{align*}
			\mathbb{E}(\rvert S_{n}^{\star}(t)-S_{n}^{\star}(v)\rvert^{p})^{1/p}\leq & \mathbb{E}(\rvert S_{n}^{\star}([t/\Delta_{n}]\Delta_{n})-S_{n}^{\star}([v/\Delta_{n}]\Delta_{n})\rvert^{p})^{1/p}\\
			& +\mathbb{E}(\rvert S_{n}^{\star}(t)-S_{n}^{\star}([t/\Delta_{n}]\Delta_{n})\rvert^{p})^{1/p}+\mathbb{E}(\rvert S_{n}^{\star}(v)-S_{n}^{\star}([v/\Delta_{n}]\Delta_{n})\rvert^{p})^{1/p}\\
			\lesssim & \Delta_{n}^{(H+1/2)\land3H}[\Delta_{n}([t/\Delta_{n}]-[v/\Delta_{n}])]{}^{1/2}+\Delta_{n}^{H+1}\\
			\leq & \Delta_{n}^{(H+1/2)\land3H}[\Delta_{n}([t/\Delta_{n}]-[v/\Delta_{n}])]{}^{1/2}+\Delta_{n}^{H+1/2}(t-v)^{1/2}\\
			\lesssim & \Delta_{n}^{(H+1/2)\land3H}(t-v)^{1/2},
		\end{align*}
		where in the last step we used that $\Delta_{n}([t/\Delta_{n}]-[v/\Delta_{n}])\lesssim t-v$
		(see e.g. Theorem 1.5.2 in \cite{BinghamGoldTeu87}).\end{proof} 
	
	\begin{remark}\label{rmk_bettererate_rough1}When $b$ is of class
		$\mathcal{C}^{4}$ with $b^{\prime}\in\mathcal{C}_{b}^{3}$ the rate
		provided in Theorem \ref{thmkey-1} can be improved. Specifically,
		under this assumption, we can show, under the same notation as above,
		that 
		\[
		\mathbb{E}(\rvert\mathcal{I}_{\mathsf{I},\mathsf{J}}^{3,3}\rvert^{p})^{1/p}\lesssim\Delta_{n}^{4H}\left[\Delta_{n}(\mathsf{J}-\mathsf{I})\right]+\Delta_{n}^{2H+1/2}\left[\Delta_{n}(\mathsf{J}-\mathsf{I})\right],
		\]
		revealing that in this situation 
		\[
		\mathbb{E}(\rvert S_{n}^{\star}(t)-S_{n}^{\star}(v)\rvert^{p})^{1/p}\leq C\Delta_{n}^{(H+1/2)\land4H}(t-v)^{p/2},\,\,\,0\leq v<t\leq T.
		\]
		For our porpouses, the estimate presented in the previous theorem
		is enough so we do not further investigate in this direction.
		
	\end{remark} 
	
	\subsection{Proof of Theorem \ref{thmerrordisGaussiannonrough}\label{subsec:Proof-of-Theoremnonrough}}
	
	In this section we will present a proof Theorem \ref{thmerrordisGaussiannonrough}.
	Our arguments heavily rely on the concept and properties of differential
	resolvents, thus we advice the reader to carefully read Subsection
	\ref{subsec:Differential-resolvents}. Recall that for $\kappa_{n}$
	as in (\ref{eq:kappa_n}), $R_{n}$ denotes its differential resolvent.
	We prepare two lemmas. The proof of the first one is a simple consequence
	of (\ref{eq:resolventequation}), Grownwall's inequality, while the
	proof of the second one further use (\ref{eq:Holdercond_X}) and the
	growth condition imposed on $b$.
	
	\begin{lemma}\label{lemmapproxresolvent}Assume that $b$ is of class
		$\mathcal{C}^{2}$ with $b^{\prime}\in\mathcal{C}_{b}^{1}$. Let $R$
		be the differential resolvent of 
		\[
		\kappa(t,A)=-\int_{[0,t]}b^{\prime}(X_{t-r})\mathbf{1}_{A}(t-r)\eta(\mathrm{d}r).
		\]
		Then for every $T>$0 there is a constant independent of $n$ such
		that 
		\begin{align}
			\sup_{0\leq s\leq v\leq t}\vert E_{n}(v,s)\vert & \leq C\sup_{0\leq s\leq T}\vert U_{s}^{n}\vert.\label{eq:R_replaceby_Rn}
		\end{align}
		
	\end{lemma}
	
	\begin{lemma}\label{estimatesincXn}Assume that $b$ is of class
		$\mathcal{C}^{2}$ with $b^{\prime}\in\mathcal{C}_{b}^{1}$ and $X^{n}$
		be as in (\ref{eq:Eulermarscheme}). Then for every $0<\lambda<H$
		and $T>0$ there is a positive random variable $\zeta_{\lambda,T}$
		only depending on $\lambda$ and $T>0$ having finite moments of all
		orders such that 
		\[
		\vert X_{t}^{n}-X_{s}^{n}\vert\leq C(\sup_{r\leq T}\vert U_{r}^{n}\vert+\vert X_{t}-X_{s}\vert)\leq C(\sup_{r\leq T}\vert U_{r}^{n}\vert+\zeta_{\lambda,T}\vert t-s\vert^{\lambda}),\,\,\,t,s\in[0,T],
		\]
		for some constant $C$ independent of $s,t$ and $n$ and $\zeta_{\lambda,T}$
		has finite moments of all order.
		
	\end{lemma}
	
	We are now ready to present our proof of Theorem \ref{thmerrordisGaussiannonrough}.
	Throughout the proof we will use the symbol $\zeta_{\lambda,T}$ to
	represent a positive random variable only depending on $\lambda$
	and $T>0$ having finite moments of all orders. Below we will also
	use the notation introduced in (\ref{eq:semilinear limit}), c.f.
	Lemma \ref{exist_unique_LinearDDE}.
	
	\begin{proof}[Proof of Theorem \ref{thmerrordisGaussiannonrough}]
		
		The proof will be esentially a consequence of the representation (\ref{eq:errordecomp2})
		and Theorem \ref{thmkey}. Our first goal is to show that we can replace
		$R_{n}$ by $R$ in relation (\ref{eq:errordecomp2}). Fix $T>0$
		and set 
		\[
		\widetilde{U}_{t}^{n}:=\int_{0}^{t}R(t,s)\mathrm{d}N_{s}^{n}=U(N^{n})_{t},\,\,\,0\leq t\leq T,
		\]
		where $N^{n}$ as in (\ref{eq:Noise_def}). By the secod part in (\ref{eq:errordecomp2}),
		Lemma \ref{lemmapproxresolvent}, and
		the properties of the differential resolvent, one easily deduces that
		\[
		\vert U_{t}^{n}-\widetilde{U}_{t}^{n}\vert\leq C\sup_{0\leq s\leq T}\vert U_{s}^{n}\vert\times\int_{0}^{T}\vert N_{s}^{n}\vert\mathrm{d}s.
		\]
		A trivial application of the the Cauchy-Schwarz inequality to the previous relation along with (\ref{eq:naiveineq})
		and Theorem \ref{thmkey}  result in 
		\begin{equation}
			\mathbb{E}(\sup_{t\leq T}\vert U_{t}^{n}-\widetilde{U}_{t}^{n}\vert^{p})\leq\Delta_{n}^{p(H+1)},\,\,p\geq1.\label{eq:firstappxnonrough}
		\end{equation}
		From this and (\ref{eq:Decomp_noise}) we conclude that
		\begin{align*}
			U_{t}^{n}= & U(N^{n})_{t}+\mathscr{O}_{p}^{u}(\Delta_{n}^{H+1})\\
			= & U(A^{n}+Y^{n})_{t}+U(D^{n})_{t}+\mathscr{O}_{p}^{u}(\Delta_{n}^{H+1})\\
			= & U(A^{n}+Y^{n})_{t}+U(D^{n}-D^{\prime,n})_{t}+\mathfrak{B}_{t}^{n}+\mathscr{O}_{p}^{u}(\Delta_{n}^{H+1}),
		\end{align*}
		where 
		\begin{align*}
			D_{t}^{\prime,n}= & \int_{0}^{t}\int_{[0,\tau]}[b(X_{s-r})-b(X_{s-\mathcal{T}(r)})]\eta(\mathrm{d}r)\mathrm{d}s\\
			= & \int_{[0,\tau]}\left[\int_{-r}^{-\mathcal{T}(r)}b(X_{s})\mathrm{d}s-\int_{t-r}^{t-\mathcal{T}(r)}b(X_{s})\mathrm{d}s\right]\eta(\mathrm{d}r).
		\end{align*}
		The first conclusion of this theorem, i.e. (\ref{eq:optimratenonrough}), now
		follows from the previous decomposition, (\ref{eq:rate_A_Bias}),
		Theorem \ref{thmkey} and the properties of the differential resolvent.
		Next, in view that 
		
		\[
		\sup_{0\leq t\leq T}\rvert D_{t}^{n}-D_{t}^{\prime,n}\rvert\leq C\Delta_{n}\sup_{0\leq s\leq T}\vert U_{s}^{n}\vert,
		\]
		we conclude as in (\ref{eq:firstappxnonrough}) that 
		\[
		U_{t}^{n}-\mathfrak{B}_{t}^{n}=U(A^{n}+Y^{n})_{t}+\mathscr{O}_{p}^{u}(\Delta_{n}^{H+1}).
		\]
		It should be clear to the reader that we are left to check that $\frac{1}{\Delta_{n}}U(A^{n}+Y^{n})$
		converges (in an appropiate way) toward $U(N)$, with the latter as
		stated in the theorem. In connection with this, our initial task is
		showing that as $n\rightarrow\infty$
		\begin{equation}
			\frac{1}{\Delta_{n}}A_{t}^{n}\overset{u.c.p}{\longrightarrow}\frac{1}{2}\int_{[0,\tau]}\int_{-r}^{(t-r)\land0}b^{\prime}(x_{0}(s))\mathrm{d}x_{0}(s)\eta(\mathrm{d}r).\label{eq:limitA}
		\end{equation}
		To see that this is indeed the case, note that by Taylor's Theorem,
		our assumptions and Lemma \ref{estimatesincXn}, we can write
		\[
		\frac{1}{\Delta_{n}}A_{t}^{n}=\frac{1}{\Delta_{n}}\int_{[0,\tau]}\int_{-\mathcal{T}(r)}^{(t-\mathcal{T}(r))\land0}\int_{\mathcal{T}(s)}^{s}b^{\prime}(x_{0}(\mathcal{T}(s))x_{0}^{\prime}(u)\mathrm{d}u\mathrm{d}s\eta(\mathrm{d}r)+\mathscr{O}_{p}^{u}(\Delta_{n}^{H}).
		\]
		Furthermore, the inner integrals of the first summand in the previous
		decomposition equal to 
		\[
		\Delta_{n}^{2}\sum_{k=-[r/\Delta_{n}]+1}^{([t/\Delta_{n}]-[r/\Delta_{n}])\land0}\int_{0}^{1}\int_{0}^{w}b^{\prime}(x_{0}(t_{k-1}))x_{0}^{\prime}(t_{k-1}+y\Delta_{n})\mathrm{d}y\mathrm{d}w+\mathrm{O}(\Delta_{n}^{2}).
		\]
		The convergence in (\ref{eq:limitA}) is inmediately obtained from this, the Dominated
		Convergence Theorem, and Dini's Theorem. 
		
		\textbf{Assume that $H>1/2$: }Theorem \ref{thmkey} and the Continuous
		Mapping Theorem yield
		\begin{equation}
			\frac{1}{\Delta_{n}}(Y^{n}+A^{n})\overset{u.c.p}{\longrightarrow}\frac{1}{2}\int_{[0,\tau]}(b(X_{\cdot-r})-b(X_{-r}))\eta(\mathrm{d}r)=N.\label{eq:conv_AnpluYn_Hgreat05}
		\end{equation}
		The boundedness of $R$ on $\Omega\times[0,T]^{2}$ and the preceeding
		relation readily indicates that
		\[
		\sup_{0\leq t\leq T}\rvert U(Y^{n}+A^{n})_{t}-U(N)_{t}\rvert\leq C\sup_{0\leq t\leq T}\rvert Y_{t}^{n}+A_{t}^{n}-N_{t}\rvert\overset{\mathbb{P}}{\rightarrow}0,
		\]
		wich is the conclusion of this theorem when $H>1/2$.
		
		\textbf{Assume that $H=1/2$: }Here we use the decomposition introduced
		in the proof of Theorem \ref{thmkey}. From (\ref{eq:first_order_Sn-1}),
		we may write 
		\[
		Y_{t}^{n}=\mathcal{Y}_{t}^{n}+\mathcal{A}_{t}^{n},
		\]
		where
		\begin{align*}
			\mathcal{A}_{t}^{n} & =\int_{[0,t]}[\mathfrak{C}_{t-r}^{n,1}+\mathfrak{C}_{t-r}^{n,2}]\eta(\mathrm{d}r),\\
			\mathcal{Y}_{t}^{n} & =\int_{[0,t]}M_{t-r}^{n}\eta(\mathrm{d}r).
		\end{align*}
		Using (\ref{eq:BM_drift1}), (\ref{eq:BM_drift2}), and (\ref{eq:limitA})
		we conclude as before that 
		\[
		\frac{1}{\Delta_{n}}U(A^{n}+\mathcal{A}^{n})\overset{u.c.p}{\longrightarrow}U(A),
		\]
		where 
		\begin{align*}
			A_{t}= & \frac{1}{2}\int_{[0,\tau]}\int_{-r}^{(t-r)\land0}b^{\prime}(x_{0}(s))\mathrm{d}x_{0}(s)\eta(\mathrm{d}r)+\frac{1}{4}\int_{[0,\tau]}\int_{0}^{t-r}b^{\prime\prime}(X_{s})\mathrm{d}s\eta(\mathrm{d}r)\\
			& +\frac{1}{2}\int_{[0,\tau]}\int_{0}^{t-r}b^{\prime}(X_{s})\int_{[0,\tau]}b(X_{s-y})\eta(\mathrm{d}y)\mathrm{d}s\eta(\mathrm{d}r).
		\end{align*}
		Considering that $N$ as in the theorem equals $A+\mathcal{Y}$ where
		\[
		\mathcal{Y}_{t}=\frac{1}{2}\int_{[0,\tau]}\int_{0}^{t-r}b^{\prime}(X_{s})\mathrm{d}W_{s}\eta(\mathrm{d}r)+\frac{1}{\sqrt{12}}\int_{[0,\tau]}\int_{0}^{t-r}b^{\prime}(X_{s})\mathrm{d}\tilde{W}_{s}\eta(\mathrm{d}r),
		\]
		and 
		\begin{equation}
			\frac{1}{\Delta_{n}}\mathcal{Y}^{n}\overset{\mathcal{F}-\mathcal{C}([0,T])}{\Longrightarrow}\mathcal{Y},\label{eq:eq:conv_AnpluYn_Brownian}
		\end{equation}
		thanks to (\ref{eq:BM_mtgpart}), in order to obtain the desired convergence,
		we need to show that 
		\[
		U^{\prime,n}:=\frac{1}{\Delta_{n}}\int_{0}^{\cdot}\mathcal{Y}_{s}^{n}\frac{\partial}{\partial s}R(\cdot,s)\mathrm{d}s\overset{\mathcal{F}-\mathcal{C}([0,T])}{\Longrightarrow}\frac{1}{\Delta_{n}}\int_{0}^{\cdot}\mathcal{Y}_{s}^{n}\frac{\partial}{\partial s}R(\cdot,s)\mathrm{d}s=:U^{\prime}.
		\]
		With this in mind, we first verify that $U^{\prime,n}$ is tight in $\mathcal{C}([0,T])$.
		For every $0\leq v\leq t\leq T$ decompose 
		\[
		U_{t}^{\prime,n}-U_{v}^{\prime,n}=\int_{v}^{t}\frac{\mathcal{Y}_{s}^{n}}{\Delta_{n}}\frac{\partial}{\partial s}R(t,s)\mathrm{d}s+\int_{0}^{v}\frac{\mathcal{Y}_{s}^{n}}{\Delta_{n}}\left[\frac{\partial}{\partial s}R(t,s)-\frac{\partial}{\partial s}R(v,s)\right]\mathrm{d}s.
		\]
		From (\ref{eq:derivative_diif_resolv}) we deduce that $\frac{\partial}{\partial s}R(\cdot,\cdot)$
		is uniformly bounded and that the second integral in the previous
		equation can be further decomposed as 
		\[
		-\int_{(0,t]}\int_{v-r}^{v\land(t-r)}\frac{\mathcal{Y}_{s}^{n}}{\Delta_{n}}b^{\prime}(X_{s})R(t,s+r)\mathrm{d}s\eta(\mathrm{d}r)-\int_{0}^{v}\frac{\mathcal{Y}_{s}^{n}}{\Delta_{n}}b^{\prime}(X_{s})\int_{[0,v-s]}[R(t,s+r)-R(v,s+r)]\eta(\mathrm{d}r)\mathrm{d}s
		\]
		These observations in conjunction with the first part of (\ref{eq:resolventequation})
		and (\ref{eq:tightness_M_BM}) result in 
		\[
		\mathbb{E}(\rvert U_{t}^{\prime,n}-U_{v}^{\prime,n}\rvert^{p})\leq C(t-v)^{p},
		\]
		as required. Therefore, we are left to show the convergence of the
		finite-dimensional distributions of $	U^{\prime,n}$. To obtain such a convergence, for
		$m\in\mathbb{N}$, set
		\begin{align*}
			U_{t}^{\prime,n,m} & :=\int_{0}^{t}\frac{\mathcal{Y}_{[s/\Delta_{m}]\Delta_{m}}^{n}}{\Delta_{n}}\frac{\partial}{\partial s}R(t,s)\mathrm{d}s,\\
			U_{t}^{\prime,m} & :=\int_{0}^{t}\mathcal{Y}_{[s/\Delta_{m}]\Delta_{m}}\frac{\partial}{\partial s}R(t,s)\mathrm{d}s.
		\end{align*}
		Note that by the properties of the stable convergence, (\ref{eq:eq:conv_AnpluYn_Brownian}),
		and the continuity of $\mathcal{Y}$ we have that: $i)$ As $n\rightarrow\infty$,
		the finite-dimensional distributions of $U^{\prime,n,m}$ converge
		$\mathcal{F}$-stably in distribution toward those of $U^{m}$; $ii)$
		$U_{t}^{m}\overset{\mathbb{P}}{\longrightarrow}U(\mathcal{Y})_{t}$
		as $m\rightarrow\infty$. Moreover, from (\ref{eq:tightness_M_BM})
		we also deduce that
		\[
		\mathbb{E}(\rvert U_{t}^{\prime,n,m}-U_{t}^{\prime,n}\rvert)\leq C\Delta_{m}^{1/2}.
		\]
		The sought-after convergence now follows as an application of Theorem
		3.21 in \cite{HauslerLuschgy15}.\end{proof} 
	We finish this section by demonstrating that \eqref{eq:limitforVn}
	is valid.
	
	\begin{proof}[Proof of \eqref{eq:limitforVn}]	Assume that $\eta$ admits a continuous density. Using this, (\ref{eq:derivative_diif_resolv}),
		and (\ref{eq:bounmomentssol}) we can decompose 
		\begin{equation}
			\begin{aligned}\mathfrak{B}_{t}^{n}= & \int_{[0,\tau]}\int_{-\mathcal{T}(r)}^{t-r}b(X_{s})[R(t,s+r)-R(t,s+\mathcal{T}(r))]\mathrm{d}s\eta(\mathrm{d}r)\\
				& -\int_{[0,\tau]}\int_{t-r}^{t-\mathcal{T}(r)}b(X_{s})R(t,s+\mathcal{T}(r))\mathrm{d}s\eta(\mathrm{d}r)\\
				& +\int_{[0,\tau]}\int_{-r}^{-\mathcal{T}(r)}b(X_{s})R(t,s+r)\mathrm{d}s\eta(\mathrm{d}r)\\
				= & \int_{0}^{t}\frac{\partial R(t,s)}{\partial s}\int_{[0,\tau]}b(X_{s-r})(r-\mathcal{T}(r))\eta(\mathrm{d}r)\mathrm{d}s\\
				& -R(t,t)\int_{[0,\tau]}b(X_{t-r})(r-\mathcal{T}(r))\eta(\mathrm{d}r)\\
				& +R(t,0)\int_{[0,\tau]}b(X_{-r})(r-\mathcal{T}(r))\eta(\mathrm{d}r)+\mathscr{O}_{p}^{u}(\Delta_{n}^{H+1}).
			\end{aligned}
			\label{eq:decomp_bias}
		\end{equation}
		Thus, by reasoning as in (\ref{eq:limitA}) and using the boundedness
		of $R$ we easily obtain that
		\begin{align}
			\frac{1}{\Delta_{n}}\mathfrak{B}_{t}^{n}\overset{u.c.p}{\longrightarrow} & -\frac{1}{2}\int_{[0,\tau]}\left[b(X_{t-r})R(t,t)-b(X_{-r})R(t,0)-\int_{0}^{t}b(X_{s-r})\frac{\partial R(t,s)}{\partial s}\mathrm{d}s\right]\eta(\mathrm{d}r).\label{eq:convBias}
		\end{align}
		Finally, by the properties of the differential resolvent, the integrand in the right-hand-side of the previous integral can be further represented as
		\[
		b(X_{t-r})-b(X_{-r})-\int_{0}^{t}[b(X_{s-r})-b(X_{-r})]\frac{\partial R(t,s)}{\partial s}\mathrm{d}s.
		\]
		The convergence in (\ref{eq:limitforVn}) is now obtained by applying this identity to the right-hand side of (\ref{eq:convBias})
		and using Fubini's Theorem. \end{proof}
	
	\subsection{Proof of Theorem \ref{themroughcase}}
	
	Here we present a proof of Theorem \ref{themroughcase}. We will use
	the following estimate 
	\begin{equation}
		\mathbb{E}(\rvert X_{t}^{n}-X_{v}^{n}-Z_{t}-Z_{v}\rvert^{p})^{1/p}\lesssim(t-v),\label{eq:incre_sol_main}
	\end{equation}
	which can be obtained by enforcing \eqref{eq:decompFLSM} and \eqref{eq:acpartZ}
	as well as Remark \ref{rmk_momentEuler} into our defition of $X^{n}$
	(see \eqref{eq:Eulermarscheme}).
	
	\begin{proof}[Proof of Theorem \ref{themroughcase} ] The main
		goal of the proof is to show that 
		\begin{equation}
			S_{n}(t)=S_{n}^{\star}(t)+\mathscr{O}_{p}^{u}(\Delta_{n}^{1\land3H}),\label{eq:approx_rough}
		\end{equation}
		where $S_{n}^{\star}$ as in (\ref{eq:def_Sn_star}). If this were
		true, we would infer from (\ref{eq:errordecomp2}), (\ref{eq:rate_A_Bias}),
		and our assumptions on $b$, that
		\[
		U_{t}^{n}=Y_{t}^{\star,n}+\int_{0}^{t}Y_{s}^{\star,n}\int_{[0,t-s]}\psi_{s}^{n}R_{n}(t,s+r)\eta(\mathrm{d}r)\mathrm{d}s+\mathscr{O}_{p}^{u}(\Delta_{n}^{1\land3H}),
		\]
		where
		\[
		Y_{t}^{\star,n}=\int_{[0,t]}S_{n}^{\star}(t-r)\eta(\mathrm{d}r).
		\]
		This along with Theorem \ref{thmkey-1} will inmediately give the
		desired result. As a first step to achieve (\ref{eq:approx_rough}),
		we derive the following non-optimal estimate 
		\begin{equation}
			\sup_{0\leq t\leq T}\mathbb{E}(\rvert S_{n}(t)\rvert^{p})^{1/p}\lesssim\Delta_{n}^{2H}.\label{eq:second_naiverate-1}
		\end{equation}
		Indeed, by arguing exactly as in the first part of the proof of Theorem
		\ref{thmkey} (see the arguments between equations (\ref{eq:first_order_Sn})
		and (\ref{eq:BV_tightness})) we have that 
		\begin{align}
			\begin{aligned}S_{n}(t)= & \int_{0}^{t}b^{\prime}(X_{\mathcal{T}(s)}^{n})[Z_{s}-Z_{\mathcal{T}(s)}]\mathrm{d}s+\frac{1}{2}\int_{0}^{t}b^{\prime\prime}(X_{\mathcal{T}(s)})[X_{s}^{n}-X_{\mathcal{T}(s)}^{n}]^{2}\mathrm{d}s+\mathscr{O}_{p}^{u}(\Delta_{n}^{3H\land1})\\
				= & L_{t}^{n}+\mathscr{O}_{p}^{u}(\Delta_{n}^{2H})
			\end{aligned}
			\label{eq:decomp1_Sn_rough}
		\end{align}
		where $L^{n}$ as in (\ref{eq:def_L_Q}). In view that
		\[
		L_{t}^{n}=\mathcal{I}_{0,[t/\Delta_{n}]}^{2}+\mathscr{O}_{p}^{u}(\Delta_{n}),
		\]
		where $\mathcal{I}_{0,[t/\Delta_{n}]}^{2}$ as in (\ref{eq:basic_decomposition_Sn}),
		we conclude from the decomposition (\ref{eq:decomp_I2}) and the estimates
		(\ref{eq:I21_rough_lpest})-(\ref{eq:norm_bound_meyet}) that 
		\begin{align*}
			\mathbb{E}(\rvert L_{t}^{n}\rvert^{p})^{1/p}\lesssim & \mathbb{E}(\rvert\mathcal{I}_{0,[t/\Delta_{n}]}^{2,3}\rvert^{p})^{1/p}+\Delta_{n}^{H+1/2}\\
			\lesssim & \sum_{i=1}^{[t/\Delta_{n}]}\rvert\int_{0}^{t_{i-1}}(t_{i-1}-u)_{+}^{\beta}\psi_{i}^{n}(u)\mathrm{d}u\rvert+\Delta_{n}^{H+1/2}\\
			\lesssim & \Delta_{n}^{2H}+\Delta_{n}^{H+1/2}
		\end{align*}
		where in the last two inequalities we further used (\ref{eq:decomp_I23})
		and Lemma (\ref{Lemma_estimates_phi}), respectively. Relation (\ref{eq:second_naiverate-1})
		now follows from this and the fact that $H<1/2$. Note this bound
		implies that 
		\begin{equation}
			\sup_{0\leq t\leq T}\mathbb{E}(\rvert U_{t}\rvert^{p})^{1/p}\lesssim\Delta_{n}^{2H}.\label{eq:second_naiverate}
		\end{equation}
		An application of the previous inequality to (\ref{eq:decomp1_Sn_rough})
		yield 
		\begin{align*}
			S_{n}(t) & =S_{n}^{\star}(t)+\frac{1}{2}E_{n}(t)+\mathscr{O}_{p}^{u}(\Delta_{n}^{3H\land1}),
		\end{align*}
		in which
		\begin{align*}
			E_{n}(t)= & \int_{0}^{t}b^{\prime\prime}(X_{\mathcal{T}(s)})[(X_{s}^{n}-X_{\mathcal{T}(s)}^{n})^{2}-(Z_{s}-Z_{\mathcal{T}(s)})^{2}]\mathrm{d}s.
		\end{align*}
		Finally, the basic relation $x^{2}-y^{2}=(x+y)(x-y)$, equations (\ref{eq:incre_sol_main})
		and (\ref{eq:bounmoment_Euler}) plus the boundedness of $b^{\prime\prime}$
		let us conclude that $E_{n}=\mathscr{O}_{p}^{u}(\Delta_{n})$, completing
		this the proof.\end{proof} 
	\subsection{Proof of Theorem \ref{thmerrordisGaussiannonrough-1}}Unlike in the previous proofs, our analysis focus directly on $N^{n}$
	and not on $S_{n}$. Recall that in this situation we are assuming that
	$\eta$ has a Lipschitz density that we will denote by $f_{\eta}$.
	
	\begin{proof}[Proof of Theorem \ref{thmerrordisGaussiannonrough-1}]
		As in the previous sections, we only need to show that under our assumptions
		\[
		\sup_{0\leq t\leq T}\mathbb{E}(\rvert N_{t}^{n}\rvert^{p})^{1/p}\lesssim\Delta_{n}^{H+1/2}.
		\]
		due to (\ref{eq:errordecomp2}). Set $H^{n}=X^{n}-B,$ and decompose
		\[
		N_{t}^{n}=L_{t}^{n,\star}+\sum_{\ell=1}^{4}A_{t}^{n,(\ell)}+\mathscr{O}_{p}^{u}(\Delta_{n}^{2H}),
		\]
		where
		\begin{align*}
			L_{t}^{n,\star} & :=\sum_{i=1}^{[t/\Delta_{n}]}\sum_{j=1}^{n}\mathbf{1}_{i\geq j+1}b^{\prime}(X_{t_{i}-t_{j}})\int_{t_{j-1}}^{t_{j}}\int_{t_{i-1}}^{t_{i}}(Z_{s-r}-Z_{t_{i}-t_{j}})\eta(\mathrm{d}r)\mathrm{d}s;\\
			A_{t}^{n,(1)} & :=\sum_{i=1}^{[t/\Delta_{n}]}\sum_{j=1}^{n}\mathbf{1}_{i\geq j+1}b^{\prime}(X_{t_{i}-t_{j}})\int_{t_{j-1}}^{t_{j}}\int_{t_{i-1}}^{t_{i}}(V_{s-r}-V_{t_{i}-t_{j}})\eta(\mathrm{d}r)\mathrm{d}s;\\
			A_{t}^{n,(2)} & :=\sum_{i=1}^{[t/\Delta_{n}]}\sum_{j=1}^{n}\mathbf{1}_{i\geq j+1}b^{\prime}(X_{t_{i}-t_{j}})\int_{t_{j-1}}^{t_{j}}\int_{t_{i-1}}^{t_{i}}(H_{s-r}^{n}-H_{t_{i}-t_{j}}^{n})\eta(\mathrm{d}r)\mathrm{d}s;\\
			A_{t}^{n,(3)} & :=\sum_{i=1}^{[t/\Delta_{n}]}\sum_{j=1}^{n}\mathbf{1}_{i\leq j-1}b^{\prime}(x_{0}(t_{i}-t_{j}))\int_{t_{j-1}}^{t_{j}}\int_{t_{i-1}}^{t_{i}}(x_{0}(s-r)-x_{0}(t_{i}-t_{j}))\eta(\mathrm{d}r)\mathrm{d}s;\\
			A_{t}^{n,(4)} & :=b^{\prime}(x_{0}(0))\sum_{i=1}^{[t/\Delta_{n}]}\int_{t_{i-1}}^{t_{i}}\int_{t_{i-1}}^{t_{i}}(X_{s-r}-X_{0})\eta(\mathrm{d}r)\mathrm{d}s.
		\end{align*}
		The leading term is $L^{n,\star}$. Let us focus on the other terms
		first. Obviously, $A_{t}^{n,(4)}=\mathscr{O}_{p}^{u}(\Delta_{n}^{H+1})$.
		Put 
		\[ A_{t}^{n,(3,+)}=\sum_{i=1}^{[t/\Delta_{n}]}\sum_{j=1}^{n}\mathbf{1}_{i\leq j-1}b^{\prime}(x_{0}(t_{i}-t_{j}))\alpha_{i,j}^{+};\,	A_{t}^{n,(3,-)}:=-\sum_{i=1}^{[t/\Delta_{n}]}\sum_{j=1}^{n}\mathbf{1}_{i\leq j-1}b^{\prime}(x_{0}(t_{i}-t_{j}))\alpha_{i,j}^{-}, \]
		where 
		\begin{align*}
			\alpha_{i,j}^{+}:=&\int_{t_{j-1}}^{t_{j}}\int_{t_{i-1}}^{t_{i}}(x_{0}(s-r)-x_{0}(t_{i}-t_{j}))\mathbf{1}_{s-r\geq t_{i}-t_{j}}\eta(\mathrm{d}r)\mathrm{d}s;\\\alpha_{i,j}^{-}:=&\int_{t_{j-1}}^{t_{j}}\int_{t_{i-1}}^{t_{i}}(x_{0}(t_{i}-t_{j})-x_{0}(s-r))\mathbf{1}_{s-r<t_{i}-t_{j}}\eta(\mathrm{d}r)\mathrm{d}s.
		\end{align*}
		We define in an analogous $A^{n,(1,\pm)}$ and $A^{n,(2,\pm)}$. Since
		$x_{0}\in C^{1}$, we can make some change of variables to deduce
		that
		\begin{equation}
			\begin{aligned}\alpha_{i,j}^{+}& =\Delta_{n}^{3}\int_{0}^{1}\int_{0}^{x}\int_{0}^{x-y}x_{0}^{\prime}(u\Delta_{n}+t_{i}-t_{j})\mathrm{d}uf_{\eta}(t_{j-1}+y\Delta_{n})\mathrm{d}y\mathrm{d}x;\\
				\alpha_{i,j}^{+-}& =\Delta_{n}^{3}\int_{0}^{1}\int_{0}^{1}\int_{0}^{x-y}x_{0}^{\prime}(t_{i}-t_{j}-u\Delta_{n})\mathrm{d}uf_{\eta}(t_{j}-y\Delta_{n})\mathrm{d}y\mathrm{d}x.
			\end{aligned}
			\label{eq:integralexp1_etadensity}
		\end{equation}
		Thus, the absolute value of the difference of these two terms is uniformly
		bounded by 
		\[
		\Delta_{n}^{3}(\sup_{\rvert u-v\rvert\leq2\Delta_{n}}\rvert x_{0}^{\prime}(u)-x_{0}^{\prime}(v)\rvert+\sup_{\rvert u-v\rvert\leq2\Delta_{n}}\rvert f_{\eta}(u)-f_{\eta}(v)\rvert)\lesssim\Delta_{n}^{4},
		\]
		due to the Lipschitz property of $x_{0}^{\prime}$ and $f_{\eta}$.
		Therefore,
		\[
		A_{t}^{n,(3)}=A_{t}^{n,(3,+)}-A_{t}^{n,(3,-)}=\mathscr{O}_{p}^{u}(\Delta_{n}^{2}).
		\]
		Next we verify that $A^{n,(2)}=\mathscr{O}_{p}^{u}(\Delta_{n}^{1+\lambda})$.
		First, we note that \eqref{eq:Eulermarscheme}, Lemma \ref{estimatesincXn},
		and Theorem \ref{thmerrordisGaussiannonrough} imply that $H_{t}^{n}=\int_{0}^{t}h_{s}^{n}\mathrm{d}s$
		with $\mathbb{E}(\sup_{t\leq T}\rvert h_{t}^{n}\rvert^{p})^{1/p}<\infty$
		as well as 
		\[
		\mathfrak{h}_{s,u}^{n,p}:=	\mathbb{E}(\rvert h_{s}^{n}-h_{u}^{n}\rvert^{p})^{1/p}\lesssim\Delta_{n}+\rvert s-u\rvert^{\lambda},\,\,1/2<\lambda<H.
		\]
		Therefore, by replacing $x_{0}^{\prime}$ by $h^{n}$ in (\ref{eq:integralexp1_etadensity})
		and applying the preceeding two properties we obtain that for all
		$1/2<\lambda<H$
		\begin{align*}
			\mathbb{E}(\sup_{0\leq t\leq T}\rvert A_{t}^{n,(2)}\rvert^{p})^{1/p} \lesssim & \Delta_{n}^{3}\sum_{i,j=1}^{[T/\Delta_{n}]\lor n}\mathbf{1}_{i\geq j+1}\int_{0}^{1}\int_{0}^{x}\int_{0}^{x-y}\mathfrak{h}_{u\Delta_{n}+t_{i}-t_{j},t_{i}-t_{j}-u\Delta_{n}}^{n,p}\mathrm{d}u\mathrm{d}y\mathrm{d}x+\Delta_{n}^{2}\\
			\lesssim&\Delta_{n}^{1+\lambda}.
		\end{align*}
		For $A^{n,(1)}$ we have instead that the density of $V$ is a Gaussian
		process satisfying that $\mathbb{E}(\vert v_{t}\vert^{2})^{1/2}\leq Ct^{H-1}$
		as well as 
		\begin{equation}
			\mathbb{E}(\vert v_{t}-v_{s}\vert^{2})^{1/2}\leq C\left\{ (t-s)^{2(\beta-1)+1}\int_{\frac{s}{t-s}}^{\infty}\{(1+y)^{\beta-1}-y^{\beta-1}\}^{2}\mathrm{d}y\right\} ^{1/2}\leq Cs^{H-1}(1\land\frac{(t-s)}{s}),\label{eq:incre_derV}
		\end{equation}
		where in the second inequality we used the fact that $(1+y)^{\beta-1}-y^{\beta-1}\sim c_{1}y^{\beta-2}$ as $y\rightarrow+\infty,$ as well as that $(1+y)^{\beta-1}-y^{\beta-1}\sim y^{\beta-1}$ when $y\rightarrow0$.Therefore, in this situation we have that 
		\[
		\mathbb{E}(\rvert v_{u\Delta_{n}+t_{i}-t_{j}}-v_{t_{i}-t_{j}-u\Delta_{n}}\rvert^{p})^{1/p}\lesssim(t_{i-1}-t_{j})^{H-2}\mathbf{1}_{i-j\geq3}\Delta_{n}u+\Delta_{n}{}^{H-1}(1-u)^{H-1}\mathbf{1}_{i-j=1,2},  \,\, u\in(0,1) .
		\]
		As a result we infer as above that 
		\[
		\mathbb{E}(\sup_{0\leq t\leq T}\rvert A_{t}^{n,(1)}\rvert^{p})^{1/p}\lesssim\Delta_{n}^{4}\sum_{j=1}^{[T/\Delta_{n}]\lor n}\sum_{i=j+2}^{[T/\Delta_{n}]\lor n}(t_{i}-t_{j})^{H-2}+\Delta_{n}^{1+H}\lesssim\Delta_{n}^{1+H},
		\]
		Since $Z$ is not absolutely continuous, the threatment of $L^{n,\star}$
		will differ slightly of the previous arguments. Recall that in this
		situation 
		\begin{equation}
			(t-u)_{+}^{\beta}=\beta\int_{u}^{t}(z-u)_{+}^{\beta-1}\mathrm{d}z,\,\,\,t,u\in\mathbb{R},\label{eq:kernel_nonrough}
		\end{equation}
		so that $Z_{t}=\beta\int_{0}^{T}\int_{u}^{t}(z-u)_{+}^{\beta-1}\mathrm{d}z\mathrm{d}W_{u}$,
		$0\leq t\leq T$. For $u\geq0$ consider 
		\begin{align*}
			\boldsymbol{I}_{i,j}^{+}(u):= & \int_{t_{j-1}}^{t_{j}}\int_{t_{i-1}}^{t_{i}}\left[(s-r-u)_{+}^{\beta}-(t_{i}-t_{j}-u)_{+}^{\beta}\right]f_{\eta}(r)\mathbf{1}_{s-r\geq t_{i}-t_{j}}\mathrm{d}r\mathrm{d}s,\\
			\boldsymbol{I}_{i,j}^{-}(u):= & \int_{t_{j-1}}^{t_{j}}\int_{t_{i-1}}^{t_{i}}\left[(t_{i}-t_{j}-u)_{+}^{\beta}-(s-r-u)_{+}^{\beta}\right]f_{\eta}(r)\mathbf{1}_{s-r<t_{i}-t_{j}}\mathrm{d}r\mathrm{d}s,
		\end{align*}
		and observe that 
		\begin{align*}
			L_{t}^{n,\star} & =\sum_{i=1}^{[t/\Delta_{n}]}\sum_{j=1}^{n}\mathbf{1}_{i\geq j+1}b^{\prime}(X_{t_{i}-t_{j}})\int_{0}^{t_{i+1}-t_{j}}\boldsymbol{I}_{i,j}^{+}(u)\mathrm{d}W_{u}-\sum_{i=1}^{[t/\Delta_{n}]}\sum_{j=1}^{n}\mathbf{1}_{i\geq j+1}b^{\prime}(X_{t_{i}-t_{j}})\int_{0}^{t_{i}-t_{j}}\boldsymbol{I}_{i,j}^{-}(u)\mathrm{d}W_{u}\\
			& =\sum_{i=2}^{[t/\Delta_{n}]}\sum_{j=1}^{n}\mathbf{1}_{i\geq j+1}b^{\prime}(X_{t_{i}-t_{j}})\int_{0}^{t_{i}-t_{j}}[\boldsymbol{I}_{i-1,j}^{+}(u)-\boldsymbol{I}_{i,j}^{-}(u)]\mathrm{d}W_{u}\\
			& +\sum_{i=2}^{[t/\Delta_{n}]+1}\sum_{j=1}^{n}\mathbf{1}_{i\geq j+1}[b^{\prime}(X_{t_{i-1}-t_{j}})-b^{\prime}(X_{t_{i}-t_{j}})]\int_{0}^{t_{i}-t_{j}}\boldsymbol{I}_{i-1,j}^{+}(u)\mathrm{d}W_{u}\\
			& +\sum_{j=1}^{n\land[t/\Delta_{n}]}b^{\prime}(X_{[t/\Delta_{n}]\Delta_{n}-t_{j}})\int_{0}^{[t/\Delta_{n}]\Delta_{n}-t_{j}}\boldsymbol{I}_{[t/\Delta_{n}],j}^{+}(u)\mathrm{d}W_{u}-b^{\prime}(X_{0})\int_{0}^{\Delta_{n}}\sum_{i=2}^{[t/\Delta_{n}]\land n}\boldsymbol{I}_{i-1,i-1}^{+}(u)\mathrm{d}W_{u}\\
			= & \sum_{m=1}^{[t/\Delta_{n}]-1}b^{\prime}(X_{t_{m}})\int_{0}^{t_{m}}\tilde{\psi}_{m}^{n}(u)\mathrm{d}W_{u}-\sum_{m=1}^{[t/\Delta_{n}]}[b^{\prime}(X_{t_{m}})-b^{\prime}(X_{t_{m-1}})]\int_{0}^{t_{m}}\varrho_{m}^{n}(u)\mathrm{d}W_{u}\\
			& +\sum_{m=([t/\Delta_{n}]-n)^{+}}^{[t/\Delta_{n}]-1}b^{\prime}(X_{t_{m}})\int_{0}^{t_{m}}\boldsymbol{I}_{[t/\Delta_{n}],[t/\Delta_{n}]-m}^{+}(u)\mathrm{d}W_{u}-b^{\prime}(X_{0})\int_{0}^{\Delta_{n}}\sum_{i=2}^{[t/\Delta_{n}]\land n}\boldsymbol{I}_{i-1,i-1}^{+}(u)\mathrm{d}W_{u}\\
			& =\sum_{\ell=1}^{4}L_{t}^{n,(\star,\ell)},
		\end{align*}
		where 
		\begin{align*}
			\tilde{\psi}_{m}^{n}(u) & :=\sum_{j=1}^{n\land([t/\Delta_{n}]-m)}[\boldsymbol{I}_{j+m-1,j}^{+}(u)-\boldsymbol{I}_{j+m,j}^{-}(u)];\\
			\varrho_{m}^{n}(u) & :=\sum_{j=1}^{n\land([t/\Delta_{n}]-m)}\boldsymbol{I}_{j+m-1,j}^{+}(u);
		\end{align*}
		The rate $\Delta_{n}^{H+1/2}$ is carried out by $L^{n,(\star,1)}$,
		so let us study the other terms first. Some simple change of variables
		and order of integration result in
		\begin{align*}
			\boldsymbol{I}_{i,j}^{+}(u) & =\beta\int_{t_{j-1}}^{t_{j}}\int_{t_{i-1}}^{t_{i}}\int_{t_{i}-t_{j}}^{s-r}(z-u)_{+}^{\beta-1}\mathrm{d}z\mathbf{1}_{s-r\geq t_{i}-t_{j}}\mathrm{d}sf_{\eta}(r)\mathrm{d}r\\
			& =\beta\int_{t_{i}-t_{j}}^{t_{i+1}-t_{j}}(z-u)_{+}^{\beta-1}\int_{t_{j-1}}^{t_{j}}\int_{t_{i-1}}^{t_{i}}f_{\eta}(r)\mathbf{1}_{z\leq s-r}\mathrm{d}sf_{\eta}(r)\mathrm{d}r\mathrm{d}z\\
			& =\beta\int_{t_{i}-t_{j}}^{t_{i+1}-t_{j}}(z-u)_{+}^{\beta-1}\int_{t_{j-1}}^{t_{j}}(t_{i}-(z+r))^{+}f_{\eta}(r)\mathrm{d}r\mathrm{d}z\\
			& =\beta\int_{t_{i}-t_{j}}^{t_{i+1}-t_{j}}(z-u)_{+}^{\beta-1}\left(\int_{0}^{t_{i+1}-t_{j}-z}xf_{\eta}(t_{i}-z-x)\mathrm{d}x\right)\mathrm{d}z.
		\end{align*}
		Similarly, 
		\[
		\boldsymbol{I}_{i,j}^{-}(u)=\beta\int_{t_{i-1}-t_{j}}^{t_{i}-t_{j}}(z-u)_{+}^{\beta-1}\left(\int_{0}^{z-(t_{i-1}-t_{j})}xf_{\eta}(x-z+t_{i-1})\mathrm{d}x\right)\mathrm{d}z.
		\]
		In particular, for all $m=1,\ldots$ and $j\leq n\land([t/\Delta_{n}]-m),$
		\begin{equation}
			\rvert\boldsymbol{I}_{j+m-1,j}^{+}(u)\rvert\leq\rVert f_{\eta}\rVert_{\infty}\Delta_{n}\rvert\psi_{m}^{n}(u)\rvert\leq\rVert f_{\eta}\rVert_{\infty}\Delta_{n}^{2}\int_{t_{m-1}}^{t_{m}}(z-u)_{+}^{\beta-1}\mathrm{d}z,\label{eq:ineI+_nonrough_smoothdens}
		\end{equation}
		as well as 
		\begin{equation}
			\rvert\boldsymbol{I}_{j+m-1,j}^{+}(u)-\boldsymbol{I}_{j+m,j}^{-}(u)-f_{\eta}(t_{j-1})\Delta_{n}\gamma_{m}^{n}(u)\rvert\leq C_{f_{\eta}}\Delta_{n}^{2}\rvert\psi_{m}^{n}(u)\rvert,\label{eq:eq:ineI+_nonrough_smoothdens_2}
		\end{equation}
		due to the Lipschitzian assumption on $f_{\eta}$ and where $\psi_{m}^{n}$
		and $\gamma_{m}^{n}(u)$ are defined as in (\ref{eq:def_rough_function-1}).
		An application of (\ref{eq:ineI+_nonrough_smoothdens}) in combination
		with (\ref{eq:estimateintegralbetadiff}) lead to
		\[
		\int_{0}^{t_{m}}\rvert\varrho_{m}^{n}(u)\rvert^{2}\mathrm{d}u\lesssim\int_{0}^{t_{m}}\rvert\psi_{m}^{n}(u)\rvert^{2}\mathrm{d}u\leq\Delta_{n}^{2(H+1)},
		\]
		reason why $L_{t}^{n,(\star,\ell)}=\mathscr{O}_{p}^{u}(\Delta_{n}^{2H})$
		for $\ell=2,4$. Now, by replacing $\psi_{i}^{n}$ by $\boldsymbol{I}_{[t/\Delta_{n}],[t/\Delta_{n}]-i}^{+}$
		in (\ref{eq:decomp_I2}) and using (\ref{eq:ineI+_nonrough_smoothdens})
		we deduce as the first part of the proof of Theorem \ref{thmkey}
		(see (\ref{eq:tighness_nonrough_I21}) and the arguments thereafter)
		that 
		\[
		\mathbb{E}(\rvert L_{t}^{n,(\star,3)}\rvert^{p})^{1/p}\lesssim\Delta_{n}^{2}.
		\]
		In the same way we obtain that 
		\[
		L_{t}^{n,(\star,1)}=\sum_{i=1}^{[t/\Delta_{n}]-1}\int_{t_{i-1}}^{t_{i}}\sum_{m=i}^{[t/\Delta_{n}]-1}b^{\prime}(X_{t_{m}})\tilde{\psi}_{m}^{n}(u)\delta W_{u}+\sum_{i=1}^{[t/\Delta_{n}]-1}\int_{t_{i-1}}^{t_{i}}\sum_{m=i}^{[t/\Delta_{n}]-1}b^{\prime\prime}(X_{t_{m}})D_{u}X_{t_{m}}\tilde{\psi}_{m}^{n}(u)\mathrm{d}u.
		\]
		Observe that (\ref{eq:eq:ineI+_nonrough_smoothdens_2}), Lemma \ref{Lemma_estimates_phi},
		and our assumption on $b$ imply that 
		\begin{equation}
			\left|\sum_{m=i}^{[t/\Delta_{n}]-1}b^{\prime}(X_{t_{m}})\tilde{\psi}_{m}^{n}(u)\right|\lesssim\Delta_{n}^{H+1/2}.\label{eq:integrandest_nonrough_smoothden}
		\end{equation}
		Hence, 
		\[
		L_{t}^{n,(\star,1)}=\sum_{i=1}^{[t/\Delta_{n}]-1}\int_{t_{i-1}}^{t_{i}}\sum_{m=i}^{[t/\Delta_{n}]-1}b^{\prime}(X_{t_{m}})\tilde{\psi}_{m}^{n}(u)\delta W_{u}+\mathscr{O}_{p}^{u}(\Delta_{n}^{H+1/2}),
		\]
		where we also made use of Lemma \ref{Malliavin_SDDE_exact}. Finally, by
		(\ref{eq:meyers_ine}), Corollary \ref{Malliavin_SDDE_exact}, and
		(\ref{eq:integrandest_nonrough_smoothden}) the $p$th moment ($p\geq1$)
		of the first term in the right-hand side of the preceeding equation is bounded up to a constant
		by 
		\[
		\left\{ \sum_{i=1}^{[t/\Delta_{n}]-1}\int_{t_{i-1}}^{t_{i}}\left(\sum_{m=i}^{[t/\Delta_{n}]-1}b^{\prime}(X_{t_{m}})\tilde{\psi}_{m}^{n}(u)\right)^{2}\mathrm{d}u\right\} ^{p/2}=\mathrm{O}(\Delta_{n}^{p(H+1/2)}),
		\]
		which is enough.\end{proof} 
	
	\subsection{Some fundamental estimates}
	
	In this section, we provide some key estimates used
	in the proof of our main results. Recall that $t_{i}=i\Delta_{n},$
	$i=0,1,\ldots$ and that $\beta:=H-1/2$. In the lemma below we will
	use the notation introduced in (\ref{eq:def_rough_function-1}). As before, the
	constants appearing below will be denoted by a generic letter $C$
	and they will be independent of $n\in\mathbb{N}$ but possibly dependent
	of $T>0$ and $H\in(0,1)$.
	
	\begin{lemma}\label{Lemma_estimates_phi}Let $T\geq t\geq0$ and
		for $t_{k-1}\leq u<t_{k},$ $k=1,\ldots,[t/\Delta_{n}]$, set 
		\[
		f_{n}(t,u):=\sum_{i=k}^{[t/\Delta_{n}]}\psi_{i}^{n}(u);\,\,g_{n}(t,u):=\sum_{i=k}^{[t/\Delta_{n}]}\rvert\gamma_{i}^{n}(u)\rvert.
		\]
		otherwise $f_{n}(t,u)=g_{n}(t,u)=h_{n}(t,u)=0$. Then, for all $u\in[0,T]$
		\begin{enumerate}
			\item If $1>H>1/2$, then
			\begin{equation}
				\rvert f_{n}(t,u)\rvert\leq C\Delta_{n};\,\,\,\rvert g_{n}(t,u)\rvert\leq C\Delta_{n}^{H+1/2},\label{eq:lemmaest_ineq1}
			\end{equation}
			and for almost all $u\in[0,T]$, $\frac{1}{\Delta_{n}}f_{n}(t,u)\rightarrow(t-u)_{+}^{\beta}.$
			\item If $\frac{1}{2}-\frac{1}{p}<H<\frac{1}{2}$, for some $p\geq1$, then
			there is a constant only depending on $T$, $p$ and $H$, such that
			\begin{equation}
				\left(\int_{0}^{t}\rvert f_{n}(t,u)\rvert^{p}\mathrm{d}u\right)^{1/p}\leq C\Delta_{n}^{H+1/2}.\label{eq:Lp_bound_fn_roguh}
			\end{equation}
			Moreover, as $n\rightarrow\infty$
			\begin{equation}
				\int_{0}^{t}\rvert f_{n}(t,u)/\Delta_{n}^{H+1/2}\rvert^{p}\mathrm{d}u\rightarrow\int_{0}^{1}\rvert f(y)\rvert^{p}\mathrm{d}y,\label{eq:limit_fn_rough}
			\end{equation}
			where 
			\[
			f(s)=\frac{s^{\beta+1}}{\beta+1}+\sum_{m\geq0}\int_{m}^{m+1}[(x+s)^{\beta}-(m+s)^{\beta}]\mathrm{d}x.
			\]
			\item Given any $\mathrm{I},\mathrm{J}=0,1,\ldots$ with $\mathrm{I}<\mathrm{J}$
			and $H\in(0,1)\backslash\{\frac{1}{2}\}$, it holds that 
			\begin{equation}
				(\int_{0}^{\mathrm{\mathrm{J}\Delta_{n}}}\rvert f_{n}(\mathrm{J}\Delta_{n},u)-f_{n}(\mathrm{I}\Delta_{n},u)\rvert^{2}\mathrm{d}u)^{1/2}\leq C\Delta_{n}^{\min\{H+1/2,1\}}[(\mathrm{J}-\mathrm{I})\Delta_{n}]^{H\lor\frac{1}{2}}.\label{eq:incr_fn_rough}
			\end{equation}
			\item Let $0<H<1/2$ and take $-1<\kappa<\frac{1}{2}-H$. Then, for every
			$\mathrm{I},\mathrm{J}=0,1,\ldots$ with $\mathrm{I}<\mathrm{J}$,
			the following estimates are attainable
			\begin{align}
				\sum_{i=\mathrm{I}}^{\mathrm{J}}\rvert\int_{0}^{t_{i}}\chi_{i}^{n}(u)\mathrm{d}u\rvert\leq & C\Delta_{n}[\Delta_{n}(\mathrm{J}-\mathrm{I})]^{2H},\label{eq:rough_cancel}\\
				\sum_{i=\mathrm{I}}^{\mathrm{J}}\rvert\int_{0}^{t_{i}}(t_{i-1}-u)_{+}^{\kappa}\psi_{i}^{n}(u)\mathrm{d}u\rvert & \leq C\Delta_{n}^{H+1/2+\kappa}\Delta_{n}(\mathrm{J}-\mathrm{I}).\label{eq:rough_int2_betakappa}
			\end{align}
		\end{enumerate}
	\end{lemma}
	
	\begin{proof}
		
		\textbf{1.} (\ref{eq:lemmaest_ineq1}) is obvious if $u\geq[t/\Delta_{n}]\Delta_{n}$,
		so take $u<[t/\Delta_{n}]\Delta_{n}$. In this situation, we can always
		write
		\begin{equation}
			\psi_{i}^{n}(u)=\beta\int_{t_{i-1}}^{t_{i}}(y-u)_{+}^{\beta-1}(t_{i}-y)\mathrm{d}y.\label{eq:phi_nonrough}
		\end{equation}
		This trivially implies that 
		\begin{equation}
			\vert f_{n}(t,u)\vert\leq\Delta_{n}(t-u)^{\beta},\,\,0\leq u\leq t.\label{eq:boundf_nbetapos}
		\end{equation}
		Now set 
		\[
		h_{n}^{(1)}(s,u,x):=\beta(\mathcal{V}(s)-u-\Delta_{n}x)_{+}^{\beta-1},\,\,0\leq s,u\leq T,0<x<1,
		\]
		in which we have let $\mathcal{V}(s)=\lceil s/\Delta_{n}\rceil\Delta_{n}$.
		A simple change of varible along with the fact that $\int_{0}^{1}(2x-1)\mathrm{d}x=0$
		give us that 
		\begin{align*}
			f_{n}(t,u) & =\Delta_{n}\int_{u}^{\mathcal{T}(t)}\int_{0}^{1}h_{n}(s,u,x)x\mathrm{d}x\mathrm{d}s+\mathfrak{I}_{1}^{n},\\
			g_{n}(t,u) & =\Delta_{n}\int_{u}^{\mathcal{T}(t)}\rvert\int_{0}^{1}[h_{n}(s,u,x)-\beta(s-u+\Delta_{n})^{\beta-1}](2x-1)\mathrm{d}x\rvert\mathrm{d}s+\rvert\mathfrak{I}_{2}^{n}\rvert,
		\end{align*}
		in which the integral is interpreted as $0$ if $u\geq\mathcal{T}(t)$
		and
		\[
		\mathfrak{I}_{1}^{n}=(u-\mathcal{T}(u))\int_{0}^{1}h_{n}(u,u,x)x\mathrm{d}x;\,\,\mathfrak{I}_{2}^{n}=(u-\mathcal{T}(u))\int_{0}^{1}h_{n}(u,u,x)(2x-1)\mathrm{d}x.
		\]
		Integration by parts yields
		\begin{align*}
			\mathfrak{I}_{1}^{n} & =\frac{(u-\mathcal{T}(u))^{\beta+2}}{(\beta+1)\Delta_{n}^{2}}=\mathrm{O}(\Delta_{n}^{\beta});\\
			\mathfrak{I}_{2}^{n} & =2\frac{(u-\mathcal{T}(u))^{\beta+2}}{(\beta+1)\Delta_{n}^{2}}-\frac{(u-\mathcal{T}(u))^{\beta+1}}{\Delta_{n}}=\mathrm{O}(\Delta_{n}^{\beta}),
		\end{align*}
		uniformly on $0\leq u\leq T$. Similarly, 
		\[
		\int_{u}^{u+\Delta_{n}}\int_{0}^{1}\rvert h_{n}(s,u,x)\rvert\mathbf{1}_{s-u\leq\Delta_{n}x}\mathrm{d}x\mathrm{d}s\leq\frac{1}{\Delta_{n}^{2}}\int_{u}^{u+\Delta_{n}}(s-u)(\mathcal{V}(s)-s)^{\beta}\mathrm{d}s\leq\Delta_{n}^{\beta}.
		\]
		Therefore, uniformly on $0\leq t,u\leq T$
		\[
		f_{n}(t,u)=\Delta_{n}\int_{u}^{\mathcal{T}(t)}\int_{0}^{1}h_{n}(s,u,x)\mathbf{1}_{s>u+\Delta_{n}x}x\mathrm{d}x\mathrm{d}s+\mathrm{O}(\Delta_{n}^{H+1/2}),
		\]
		as well as 
		\begin{equation}
			\rvert g_{n}(t,u)\rvert\leq3\Delta_{n}\int_{0}^{1}\int_{u}^{t}\rvert h_{n}(s,u,x)-(s-u+\Delta_{n})_{+}^{\beta-1}\rvert\mathbf{1}_{s>u+\Delta_{n}x}\mathrm{d}s\mathrm{d}x+C\Delta_{n}^{H+1/2},\label{eq:inegn_lemma}
		\end{equation}
		where we also used that $0<\beta<1 $and that $\mathcal{V}(s)-u-\Delta_{n}x\leq s-u+\Delta_{n}$. Now, since for all $0<x<1$ and $s>u$
		\begin{align*}
			\beta(s-u+\Delta_{n}(1-x))^{\beta-1}\mathbf{1}_{s>u+\Delta_{n}x}\leq & h_{n}(s,u,x)\mathbf{1}_{s>u+\Delta_{n}x}\leq\beta(s-u-\Delta_{n}x)^{\beta-1}\mathbf{1}_{s>u+\Delta_{n}x};
		\end{align*}
		then for fixed $u<t$ and $n$ enough large, we get that 
		\[
		\Delta_{n}\int_{0}^{1}(t-u+\Delta_{n}(1-x))^{\beta}\mathrm{d}s+\mathrm{O}(\Delta_{n}^{H+1/2})\leq f_{n}(t,u)\leq\Delta_{n}\int_{0}^{1}(t-u-\Delta_{n}x)^{\beta}x\mathrm{d}x\mathrm{d}s+\mathrm{O}(\Delta_{n}^{H+1/2}).
		\]
		A simple application of the Dominated Convergence Theorem gives that
		$\frac{1}{\Delta_{n}}f_{n}(t,u)\rightarrow(t-u)^{\beta}$ as claimed.
		Analogously, the integral on the right-hand side of \eqref{eq:inegn_lemma} can be estimated from above by 
		\[ \beta\int_{0}^{1}\int_{u+\Delta_{n}x}^{t+\Delta_{n}}[(s-u-\Delta_{n}x)^{\beta-1}-(s-u+\Delta_{n}x)^{\beta-1}]\mathrm{d}s\mathrm{d}x, \]
		which is in turn bounded from above by $2\int_{0}^{1}(2x)^{\beta}\mathrm{d}x\Delta_{n}^{\beta}$. Plugging these estimates into
		(\ref{eq:inegn_lemma}) concludes the proof of this part. 
		
		\textbf{2. }Fix $p\geq1$, such that $\frac{1}{2}-\frac{1}{p}<H<\frac{1}{2}$. We
		start by making the change of variables $y=(t_{k}-u)/\Delta_{n}$
		to get that 
		\[
		\int_{0}^{t}\rvert f_{n}(t,u)\rvert^{p}\mathrm{d}u=\Delta_{n}\sum_{k=1}^{[t/\Delta_{n}]}\int_{0}^{1}\left|\frac{y^{H+1/2}}{H+1/2}+\sum_{i=k+1}^{[t/\Delta_{n}]}\psi_{i}^{n}(t_{k}-y\Delta_{n})\right|^{p}\mathrm{d}y.
		\]
		Since for $i\geq k$
		\[
		\psi_{i}^{n}(t_{k}-y\Delta_{n})=\Delta_{n}^{H+1/2}\int_{i-1-k}^{i-k}\left[(x+y)^{\beta}-(i-1-k+y)^{\beta}\right]\mathrm{d}x,
		\]
		we further get that 
		\[
		\int_{0}^{t}\rvert f_{n}(t,u)\rvert^{p}\mathrm{d}u=\Delta_{n}^{p(H+1/2)}\Delta_{n}\sum_{k=1}^{[t/\Delta_{n}]}\int_{0}^{1}\rvert w_{k}^{n}(t,y)\rvert^{p}\mathrm{d}y
		\]
		\[
		w_{k}^{n}(t,y)=\frac{y^{H+1/2}}{H+1/2}+\sum_{m=0}^{[t/\Delta_{n}]-(k+1)}\int_{m}^{m+1}(x+y)^{\beta}-(m+y)^{\beta}\mathrm{d}x.
		\]
		Note now that in view that $\beta<0$ and $p\beta>-1$, it holds that
		uniformly on $k$
		\begin{equation}
			\int_{0}^{1}\rvert w_{k}^{n}(t,y)\rvert^{p}\mathrm{d}y\leq C\int_{0}^{1}\rvert y^{\beta+1}+y^{\beta}+\sum_{m\geq1}m{}^{\beta-1}\rvert^{p}\mathrm{d}y\leq C\int_{0}^{1}\rvert1+y^{\beta}\rvert^{p}\mathrm{d}y<\infty.\label{eq:boundh_n}
		\end{equation}
		This relation trivially implies (\ref{eq:Lp_bound_fn_roguh}). Now
		observe that 
		\[
		\int_{0}^{t}\rvert f_{n}(t,u)/\Delta_{n}^{H+1/2}\rvert^{p}\mathrm{d}u=\int_{0}^{[t/\Delta_{n}]\Delta_{n}}\int_{0}^{1}\rvert w_{[u/\Delta_{n}]}^{n}(t,y)\rvert^{p}\mathrm{d}y\mathrm{d}u.
		\]
		and for fixed $u<t$, $w_{[u/\Delta_{n}]}^{n}(t,y)\rightarrow f(y)$. Hence, (\ref{eq:limit_fn_rough}) is achieved by this,
		(\ref{eq:boundh_n}) and the Dominated Convergence Theorem.
		
		\textbf{3. }Suppose first that $H>1/2$ and set $I_{\mathrm{I},\mathrm{J}}:=\int_{0}^{\mathrm{J}\Delta_{n}}\vert f_{n}(\mathrm{J}\Delta_{n},u)-f_{n}(\mathrm{I}\Delta_{n},u)\vert^{2}\mathrm{d}u$.
		In the light of (\ref{eq:boundf_nbetapos}) and (\ref{eq:phi_nonrough})
		we get in this situation that 
		\begin{align*}
			I_{\mathrm{I},\mathrm{J}} & \leq\Delta_{n}^{2}[(\mathrm{J}-\mathrm{I})\Delta_{n}]^{2\beta+1}\\
			& +\Delta_{n}^{2}\int_{0}^{\mathrm{I}\Delta_{n}}\vert(\mathrm{J}\Delta_{n}-u)^{\beta}-(\mathrm{I}\Delta_{n}-u)^{\beta}\vert^{2}\mathrm{d}u\\
			& \leq C\Delta_{n}^{2}[(\mathrm{J}-\mathrm{I})\Delta_{n}]^{2H},
		\end{align*}
		where in the last inequality we further used (\ref{eq:estimateintegralbetadiff}).
		Now suppose that $H<1/2$. Using that for $0<u<\mathrm{I}$ it holds
		that
		\begin{align*}
			f_{n}(\mathrm{J}\Delta_{n},u\Delta_{n})-f_{n}(\mathrm{I}\Delta_{n},u\Delta_{n}) & =\Delta_{n}^{H+1/2}\sum_{m=\mathrm{I}+1}^{\mathrm{J}}\int_{m}^{m+1}\{(x-u)^{\beta}-(m-u)^{\beta}\}\mathrm{d}x.\\
			& =:\Delta_{n}^{H+1/2}\sum_{m=\mathrm{I}+1}^{\mathrm{J}}c_{m}(u),
		\end{align*}
		with
		\begin{equation}
			\vert c_{m}(u)\vert\leq(m-u)^{\beta}-(m+1-u)^{\beta},\label{eq:boundam}
		\end{equation}
		along with (\ref{eq:lemmaest_ineq1}) and (\ref{eq:estimateintegralbetadiff}),
		we obtain the following inequalities
		\begin{align*}
			I_{\mathrm{I},\mathrm{J}} & \leq C\Delta_{n}^{2H+1}\left((\mathrm{J}-\mathrm{I})\Delta_{n}+\Delta_{n}\int_{0}^{\mathrm{I}}\vert(\mathrm{J}-u)^{\beta}-(\mathrm{I}-u)^{\beta}\vert^{2}\mathrm{d}u\right)\\
			& \leq C\Delta_{n}^{2H+1}\left((\mathrm{J}-\mathrm{I})\Delta_{n}+\Delta_{n}(\mathrm{J}-\mathrm{I}){}^{2H}\right)\\
			& \leq C\Delta_{n}^{2H+1}(\mathrm{J}-\mathrm{I})\Delta_{n},
		\end{align*}
		where in the last step we used the fact that $\mathrm{J}-\mathrm{I}\geq1$
		and $2H-1<0$. This concludes the proof of (\ref{eq:incr_fn_rough}).
		
		\textbf{4. }As in 2. we may write 
		\[
		\int_{0}^{t_{i}}(t_{i-1}-u)_{+}^{\kappa}\psi_{i}^{n}(u)\mathrm{d}u=\Delta_{n}^{H+\frac{1}{2}+\kappa+1}\sum_{m=1}^{i-1}\int_{0}^{1}a_{m}^{\kappa,\beta}(y)\mathrm{d}y,
		\]
		and 
		\[
		\frac{1}{\Delta_{n}^{2H+1}}\int_{0}^{t_{i}}\chi_{i}^{n}(u)\mathrm{d}u=\frac{1}{2(2\beta+1)(2\beta+2)}+\sum_{m=1}^{i-1}\int_{0}^{1}\left[a_{m}^{\beta,\beta}(y)+\frac{1}{2}b_{m}(y)\right]\mathrm{d}y,
		\]
		where 
		\begin{align*}
			a_{m}^{\kappa,\beta}(y) & :=(m+y-1)^{\kappa}\int_{m-1}^{m}\left[(x+y)^{\beta}-(m+y-1)^{\beta}\right]\mathrm{d}x\\
			b_{m}(y) & :=\int_{m-1}^{m}\left[(x+y)^{\beta}-(m+y-1)^{\beta}\right]^{2}\mathrm{d}x.
		\end{align*}
		By using that 
		\[
		\rvert a_{m}^{\kappa,\beta}(y)\rvert\leq C(y^{\kappa}(1+y)^{\beta+1}\mathbf{1}_{m=1}+(m-1)^{\kappa+\beta-1}\mathbf{1}_{m\geq2})
		\]
		we concluce that $\sum_{m\geq1}\int_{0}^{1}\rvert a_{m}^{\kappa,\beta}(y)\rvert\mathrm{d}y<\infty$
		which easily implies (\ref{eq:rough_int2_betakappa}). Now to see
		that the bound (\ref{eq:rough_cancel}) holds, observe that for all
		$m=1,2,\ldots$
		\begin{align*}
			\int_{0}^{1}\left[a_{m}^{\beta,\beta}(y)+\frac{1}{2}b_{m}(y)\right]\mathrm{d}y & =\frac{1}{2}\int_{0}^{1}\int_{m-1}^{m}\left[(x+y)^{2\beta}-(m+y-1)^{2\beta}\right]\mathrm{d}x\mathrm{d}y\\
			& =\frac{1}{2(2\beta+1)}\int_{m-1}^{m}\left[(x+1)^{2\beta+1}-x{}^{2\beta+1}\right]\mathrm{d}x\\
			& -\frac{1}{2(2\beta+1)}\left[m^{2\beta+1}-(m-1)^{2\beta+1}\right],
		\end{align*}
		so that 
		\begin{align*}
			\frac{1}{\Delta_{n}^{2H+1}}\int_{0}^{t_{i}}\chi_{i}^{n}(u)\mathrm{d}u & =\frac{1}{2(2\beta+1)(2\beta+2)}+\frac{1}{2(2\beta+1)}\left[\int_{i-1}^{i}x^{2\beta+1}\mathrm{d}x-\int_{0}^{1}x^{2\beta+1}\mathrm{d}x\right]\\
			& -\frac{1}{2(2\beta+1)}(i-1)^{2\beta+1}\\
			& =\frac{1}{2(2\beta+1)}\int_{i-1}^{i}\left[x^{2\beta+1}-(i-1)^{2\beta+1}\right]\mathrm{d}x\\
			& \leq C_{H}(\int_{0}^{1}x^{2\beta+1}\mathrm{d}x+\int_{i-1}^{i}\left[x^{2\beta+1}-(x-1)^{2\beta+1}\right]\mathrm{d}x\mathbf{1}_{i\geq1}),
		\end{align*}
		because $2\beta+1=2H>0$. The desired inequality now follows easily
		from this estimate.
		
	\end{proof} 
	
	\bibliographystyle{plain}
	\bibliography{bibSept24}

\begin{thebibliography}{10}

\bibitem{ArrioHuMohPap07}
M.~Arriojas, Y.~Hu, S-E Mohammed, and G.~Pap.
\newblock A delayed black and scholes formula.
\newblock {\em Stochastic Analysis and Applications}, 25(2):471--492, 2007.

\bibitem{AzmooSotViitYa14}
Ehsan Azmoodeh, Tommi Sottinen, Lauri Viitasaari, and Adil Yazigi.
\newblock Necessary and sufficient conditions for h\"{o}lder continuity of
  gaussian processes.
\newblock {\em Statistics \& Probability Letters}, 94:230--235, 2014.

\bibitem{BaoYinYuan16}
J.~Bao, G.~Yin, and C.~Yuan.
\newblock {\em Asymptotic Analysis for Functional Stochastic Differential
  Equations}.
\newblock SpringerBriefs in Mathematics. Springer International Publishing,
  2016.

\bibitem{BasseBN11}
O.~E. Barndorff-Nielsen and A.~Basse-O'Connor.
\newblock Quasi {O}rnstein-{U}hlenbeck processes.
\newblock {\em Bernoulli}, 17:916--941, 2011.

\bibitem{BayerfrizFukaGatherakJacqRosem23}
Christian Bayer, Peter~K. Friz, Masaaki Fukasawa, Jim Gatheral, Antoine
  Jacquier, and Mathieu Rosenbaum.
\newblock {\em Rough Volatility}.
\newblock Society for Industrial and Applied Mathematics, Philadelphia, PA,
  2023.

\bibitem{BenLundPakk16}
M.~Bennedsen, A.~Lunde, and M.~S. Pakkanen.
\newblock Decoupling the short- and long-term behavior of stochastic
  volatility.
\newblock {\em Journal of Financial Econometrics}, 01 2021.

\bibitem{BesaluRov12}
Mireia Besal\'u and Carles Rovira.
\newblock Stochastic delay equations with non-negativity constraints driven by
  fractional brownian motion.
\newblock {\em Bernoulli}, 18(1):24--45, 2012.

\bibitem{BinghamGoldTeu87}
N.~H. Bingham, C.~M. Goldie, and J.~L. Teugels.
\newblock {\em Regular Variation}.
\newblock Encyclopedia of Mathematics and its Applications. Cambridge
  University Press, 1987.

\bibitem{BOUFHajji11}
B.~Boufoussi and S.~Hajji.
\newblock Functional differential equations driven by a fractional brownian
  motion.
\newblock {\em Computers and Mathematics with Applications}, 62(2):746--754,
  2011.

\bibitem{BUCKWAR00}
Evelyn Buckwar.
\newblock Introduction to the numerical analysis of stochastic delay
  differential equations.
\newblock {\em Journal of Computational and Applied Mathematics},
  125(1):297--307, 2000.
\newblock Numerical Analysis 2000. Vol. VI: Ordinary Differential Equations and
  Integral Equations.

\bibitem{BurdzySwanson10}
K.~Burdzy and J.~Swanson.
\newblock A change of variable formula with it{\^o} correction term.
\newblock {\em The Annals of Probability}, 38(5):1817--1869, 2010.

\bibitem{CaraballoGarridoTanigu11}
T.~Caraballo, M.J. Garrido-Atienza, and T.~Taniguchi.
\newblock The existence and exponential behavior of solutions to stochastic
  delay evolution equations with a fractional brownian motion.
\newblock {\em Nonlinear Analysis: Theory, Methods \& Applications},
  74(11):3671--3684, 2011.

\bibitem{DiNunnoKubMishuYour23}
G.a Di~Nunno, K.~Kubilius, Y.~Mishura, and A.~Yurchenko-Tytarenko.
\newblock From constant to rough: A survey of continuous volatility modeling.
\newblock {\em Mathematics}, 11(19), 2023.

\bibitem{DonKutz96}
P.~Donnelly and T.~G. Kurtz.
\newblock A countable representation of the fleming-viot measure-valued
  diffusion.
\newblock {\em The Annals of Probability}, 24(2):698--742, 1996.

\bibitem{FerrantRovira06}
M.~Ferrante and C.~Rovira.
\newblock Stochastic delay differential equations driven by fractional brownian
  motion with hurst parameter h>1/2.
\newblock {\em Bernoulli}, 12(1):85--100, 2006.

\bibitem{FerrantRovira10}
M.~Ferrante and C.~Rovira.
\newblock Convergence of delay differential equations driven by fractional
  brownian motion.
\newblock {\em Journal of Evolution Equations}, 10(4):761--783, 2010.

\bibitem{GarzonTindelTorres19}
J.a Garz{\'o}n, S.y Tindel, and S.~Torres.
\newblock Euler scheme for fractional delay stochastic differential equations
  by rough paths techniques.
\newblock {\em Acta Mathematica Scientia}, 39(3):747--763, 2019.

\bibitem{GathJaiRos18}
J.~Gatheral, T.~Jaisson, and M.~Rosenbaum.
\newblock Volatility is rough.
\newblock {\em Quantitative Finance}, 18(6):933--949, 2018.

\bibitem{GradNourdin09}
M.~Gradinaru and I.~Nourdin.
\newblock {Milstein's type schemes for fractional SDEs}.
\newblock {\em {Annales de l'Institut Henri Poincar\'e, Probabilit\'es et
  Statistiques}}, 45(4):1085 -- 1098, 2009.

\bibitem{GripenbergLondenStaffans90}
G.~Gripenberg, S.~O. Londen, and O.~Staffans.
\newblock {\em Volterra Integral and Functional Equations}.
\newblock Encyclopedia of Mathematics and its Applications. Cambridge
  University Press, 1990.

\bibitem{HauslerLuschgy15}
Erich H{\"{a}}usler and Harald Luschgy.
\newblock {\em Stable Convergence and Stable Limit Theorems}.
\newblock Probability Theory and Stochastic Modelling. Springer International
  Publishing Switzerland, 2015.

\bibitem{HoffmanMuller06}
N.~Hofmann and T.~M{\"u}ller-Gronbach.
\newblock A modified milstein scheme for approximation of stochastic delay
  differential equations with constant time lag.
\newblock {\em Journal of Computational and Applied Mathematics},
  197(1):89--121, 2006.

\bibitem{HuMohaYan04}
Y.. Hu, S.~A. Mohammed, and Feng Yan.
\newblock {Discrete-time approximations of stochastic delay equations: The
  Milstein scheme}.
\newblock {\em The Annals of Probability}, 32(1A):265 -- 314, 2004.

\bibitem{HuLiuNualart16}
Yaozhong Hu, Yanghui Liu, and David Nualart.
\newblock {Rate of convergence and asymptotic error distribution of Euler
  approximation schemes for fractional diffusions}.
\newblock {\em The Annals of Applied Probability}, 26(2):1147 -- 1207, 2016.

\bibitem{JacPrott98}
J.~Jacod and P.~Protter.
\newblock {Asymptotic error distributions for the Euler method for stochastic
  differential equations}.
\newblock {\em The Annals of Probability}, 26(1):267 -- 307, 1998.

\bibitem{JacProt11}
J.~Jacod and P.~Protter.
\newblock {\em Discretization of Processes}.
\newblock Stochastic Modelling and Applied Probability. Springer Berlin
  Heidelberg, 2011.

\bibitem{JacShri02}
J.~Jacod and A.N. Shiryaev.
\newblock {\em Limit Theorems for Stochastic Processes}.
\newblock Springer Berlin Heidelberg, 2002.

\bibitem{Kallenberg02}
O.~Kallenberg.
\newblock {\em Foundations of Modern Probability}.
\newblock Probability and Its Applications. Springer New York, 2002.

\bibitem{KaratzasShreve91}
I.~Karatzas and S.E. Shreve.
\newblock {\em Brownian Motion and Stochastic Calculus}.
\newblock Graduate Texts in Mathematics. Springer, 2nd. edition, 1991.

\bibitem{YuriSwWu05}
Yuriy Kazmerchuk, Anatoliy Swishchuk, and Jianhong Wu.
\newblock A continuous-time garch model for stochastic volatility with delay.
\newblock {\em Canadian Applied Mathematics Quarterly}, 13(2):123--149, 2005.

\bibitem{YuriSwWu07}
Yuriy Kazmerchuk, Anatoliy Swishchuk, and Jianhong Wu.
\newblock The pricing of options for securities markets with delayed response.
\newblock {\em Mathematics and Computers in Simulation}, 75(3):69--79, 2007.

\bibitem{KloedenShardlow12}
P.~E. Kloeden and T.~Shardlow.
\newblock The milstein scheme for stochastic delay differential equations
  without using anticipative calculus.
\newblock {\em Stochastic Analysis and Applications}, 30(2):181--202, 2012.

\bibitem{KuchlerPlaten07}
U.~K{\"u}chler and E.~Platen.
\newblock Time delay and noise explaining cyclical fluctuations in prices of
  commodities.
\newblock 2007.

\bibitem{KurtzPotter91}
T.~G. Kurtz and P.~Protter.
\newblock Wong-zakai corrections, random evolutions, and simulation schemes for
  sde's.
\newblock In Eddy Mayer-Wolf, Ely Merzbach, and Adam Shwartz, editors, {\em
  Stochastic Analysis}, pages 331--346. Academic Press, 1991.

\bibitem{LeonTindel12}
J.~A. Le{\'o}n and S.~Tindel.
\newblock Malliavin calculus for fractional delay equations.
\newblock {\em Journal of Theoretical Probability}, 25(3):854--889, 2012.

\bibitem{LiuTindel19}
Y.~Liu and S.~Tindel.
\newblock {First-order Euler scheme for SDEs driven by fractional Brownian
  motions: The rough case}.
\newblock {\em The Annals of Applied Probability}, 29(2):758 -- 826, 2019.

\bibitem{MAHMOUDITahmaseb22}
F.~Mahmoudi and M.~Tahmasebi.
\newblock The convergence of a numerical scheme for additive fractional
  stochastic delay equations with $h>1/2$.
\newblock {\em Mathematics and Computers in Simulation}, 191:219--231, 2022.

\bibitem{NeueNourdin07}
A.~Neuenkirch and I.~Nourdin.
\newblock Exact rate of convergence of some approximation schemes associated to
  sdes driven by a fractional {B}rownian motion.
\newblock {\em Journal of Theoretical Probability}, 20(4):871--899, 2007.

\bibitem{Neuenkirch06}
Andreas Neuenkirch.
\newblock Optimal approximation of sde's with additive fractional noise.
\newblock {\em Journal of Complexity}, 22(4):459--474, 2006.

\bibitem{NovikovValkeila99}
Alexander Novikov and Esko Valkeila.
\newblock On some maximal inequalities for fractional brownian motions.
\newblock {\em Statistics \& Probability Letters}, 44(1):47--54, 1999.

\bibitem{Nualart06}
David Nualart.
\newblock {\em {The {M}alliavin Calculus and Related Topics}}.
\newblock Probability and its Applications. Springer Berlin Heidelberg, 2nd
  edition, 2006.

\bibitem{Rihan2021}
F.A. Rihan.
\newblock {\em Delay Differential Equations and Applications to Biology}.
\newblock Springer Nature Singapore, 2021.

\bibitem{Sauri25}
O.~Sauri.

\bibitem{Zale98}
M.~Z{\"a}hle.
\newblock Integration with respect to fractal functions and stochastic
  calculus. i.
\newblock {\em Probability Theory and Related Fields}, 111(3):333--374, 1998.

\bibitem{ZhangGandHu09}
H.~Zhang, S.~Gan, and L.~Hu.
\newblock The split-step backward euler method for linear stochastic delay
  differential equations.
\newblock {\em Journal of Computational and Applied Mathematics},
  225(2):558--568, 2009.

\bibitem{ZhangYuan21}
S.~Zhang and C.~Yuan.
\newblock Stochastic differential equations driven by fractional brownian
  motion with locally lipschitz drift and their implicit euler approximation.
\newblock {\em {Proceedings of the Royal Society of Edinburgh: Section A
  Mathematics}}, 151(4):1278--1304, 2021.

\end{thebibliography}
	
\end{document}